\documentclass[12pt]{article}
\usepackage{graphicx}
\usepackage{hyperref}
\usepackage{amsmath,amssymb,amsthm}
\usepackage{enumerate}

\numberwithin{equation}{section}
\newtheorem{theorem}{Theorem}
\newtheorem{proposition}[theorem]{Proposition}
\newtheorem{lemma}[theorem]{Lemma}
\newtheorem{coro}[theorem]{Corollary}

\newtheorem{remark}[theorem]{Remark}

\newcommand{\R}{{\mathbf R}}
\newcommand{\Z}{{\mathbf Z}}
\newcommand{\N}{{\mathbf N}}
\newcommand{\Ex}{\mathbf{E} }
\newcommand{\F}{\mathcal{F} }
\newcommand{\E}{\mathcal{E} }

\title{Generalizations of Ho-Lee's binomial interest rate 
model I: from one- to multi-factor}

\author{Jir\^o Akahori$^1$,  Hiroki Aoki$^2$, 
and Yoshihiko Nagata$^3$ \\[12pt]
$^1$Department of Mathematical Sciences, \\
Ritsumeikan University \\
1-1-1 Nojihigashi, Kusatsu, Shiga 525-8577, Japan \\
e-mail: akahori@se.ritsumei.ac.jp \\[4pt]
$^2$Department of Mathematics, \\
Tokyo University of Science \\ 
Noda, Chiba, 278-8510, Japan \\
e-mail: aoki\_hiroki@ma.noda.tus.ac.jp \\[4pt]
$^3$Risk Management Department, \\
Mizuho Trust \& Banking Co. Ltd. \\
1-2-1 Yaesu, Tokyo, 103-8670, Japan \\
e-mail: yoshihiko\_nagata@tz.mizuho-tb.co.jp}

\date{This version: June 2006}

\begin{document}
\maketitle

\begin{abstract}
In this paper a multi-factor generalization of Ho-Lee 
model is proposed. 
In sharp contrast to the classical Ho-Lee, 
this generalization allows for
those movements other than parallel shifts, 
while it still is described by a recombining tree,
and is stationary to be compatible 
with principal component analysis.
Based on the model, generalizations of 
duration-based hedging are proposed. 
A continuous-time limit of the model is also discussed.
\end{abstract}
\tableofcontents

\ 

\noindent\thanks{{\bf Keywords}: Ho-Lee model, duration, 
multi-factor, recombining tree, 
stationarity, forward rate, drift condition}

\noindent\thanks{{\bf Mathematical Subject Classification 2000}
:91B28 60G50}

\noindent\thanks{{\bf Journal of Economic Literature 
Classification System}:G12}

\noindent\thanks{This research was supported by Open Research Center Project for Private Universities: matching fund subsidy from MEXT, 
2004-2008
and also by Grants-in-Aids for Scientific Research (No. 18540146)
from the Japan Society for Promotion of Sciences.}

\section{Introduction}\label{intro1}
\subsection{The aim of the present paper/missing rings between
Ho-Lee and HJM}
It is almost thirty years since 
Vasicek's paper on term structure of interest rates \cite{Vas}
was published, and is almost twenty years 
since the paper on binary interest rate by Ho and Lee \cite{HL}
appeared. 
Most of the papers on interest rate modeling in these 20-30 years 
were in anyway 
coming from either Vasicek or Ho-Lee. 

Vasicek initiated so called {\em spot rate models},
followed by Cox, Ingersoll and Ross \cite{CIR} and 
Hull and White \cite{H-W} among others, and
its multi-factor generalizations were
proposed, for example,
by Longstaff and Schwartz \cite{Lo-Sch}, 
and by Duffie and Kan \cite{D-K} they 
are unified as Affine Term Structure Models (ATSM).
Recently modifications of ATSMs 
like Quadratic Term Structure Models (QTSM)
are extensively studied (see, e.g., 
Chen, Filipovi\'c, and Poor \cite{CFP} and the references therein).
In a word, they are {\em alive}. 

After Ho-Lee, on the other hand, 
there was a major progress made by 
Heath, Jarrow and Morton \cite{HJM1} and \cite{HJM2};
the shift to 
forward rate models in a continuous-time framework\footnote{Ho and Lee themselves proposed a multi-factor generalization 
within a recombining-tree framework \cite{HL2}, 
but actually it is not that much so far.}.
The success of HJM was so great that  
Ho-Lee model almost lost its position 
as a teacher of HJM. 
In many textbooks Ho-Lee model is treated as an
almost trivial example of HJM or even worse, 
spot rate models. 

Starting from a careful study of the original Ho-Lee,
the paper will shed lights on 
what are missed and lost 
in the way between Ho-Lee and HJM or 
after HJM. 
We will liberate Ho-Lee model 
from the rule of 
parallel shifts
by shifting to forward rate modeling.
In contrast to HJM, we preserve 
the traditions of recombining tree and stationarity, 
and hence consistency with principal component analysis.

As for applications, we will concentrate on 
generalizations of the duration-based hedging.  

Most of our results presented here could be already well-known 
among practitioners, but in anyway 
they have not been explicitly stated.
We do this 
to recover missing rings between practitioners and academics.
Though the style of the presentation is 
that of mathematicians, we have carefully avoided 
too abstract mathematics. 

\subsection{Sensitivity analysis of 
interest rate risks
/What do we mean by ``multi-factor"?}\label{sensitivity}
\begin{quote}
The key assumption underlying the duration-based hedging scheme 
is that all interest rate change by the same amount.
This means that only parallel shifts in the term structure are allowed for.
(J. Hull \cite{Hull}) 
\end{quote}
The present value of a stream of cash flow is 
a function of the current term structure of interest rate. 
The {\em duration} and the {\em convexity} are the most classical
sensitivity criteria to handle the interest rate risks 
(see e.g. \cite[section 4.8 \& 4.9]{Hull}). 

Let us review the context. 
Suppose that, for $ T= T_1, \ldots, T_k $ in the future, 
we will have $ \mathbf{CF} (T) $ of deterministic cash flow, while 
the term structure of interest rates are given 
as $ \mathbf{r}(T) $. 
Then the present value at time $ t $ ($<T_1$) of the cash flow is 
\begin{equation}\label{cashflow}
\mathrm{PV}_t \equiv 
\mathrm{PV} (\mathbf{r} (T_1), \cdots, \mathbf{r} (T_k) ; t ) =
\sum_{l=1}^k \mathbf{CF} (T_l) e^{-\mathbf{r}(T_l) (T_l-t)}.
\end{equation}
The {\em duration} in the usual sense 
\begin{equation*}
(\text{duration}) \equiv (\text{dur}):= \frac{1}{\mathrm{PV}} 
\sum_{l=1}^k (T_l -t)\mathbf{CF} (T_l) e^{-\mathbf{r}(T_l) (T_l-t)}
\end{equation*}
can be obtained as 
a derivative of the multivariate function $ \log \mathrm{PV} $ 
in the direction of $ {\mathbf{ps}} \equiv (1, \ldots, 1)$:
\begin{equation}\label{der1}
\begin{split}
D_{\mathbf{ps}} :=- \partial_{{\mathbf{ps}}} (\log \mathrm{PV}) 
&= - \frac{1}{\mathrm{PV}}
\lim_{\epsilon \to 0} \epsilon^{-1} 
\{ \mathrm{PV} ( \mathbf{r} + \epsilon\, \mathbf{ps} ;t) 
- \mathrm{PV}(\mathbf{r};t)\} = (\text{dur}). 
\end{split}
\end{equation}
The {\em convexity} also comes from the second order derivative 
\begin{equation}\label{der2}
D^2_{\mathbf{ps}} := 
\frac{1}{\mathrm{PV}}  \partial^2_{{\mathbf{ps}}} (\mathrm{PV})
= - \frac{1}{\mathrm{PV}} 
\sum_{l=1}^k (T_l -t)^2\mathbf{CF} (T_l) e^{-\mathbf{r}(T_l) (T_l-t)}. 
\end{equation}

Suppose that the random dynamics 
of the interest rates 
over the period $ [t, t + \Delta t ] $
is described by a real valued 
random variable $ \Delta w_t $:
\begin{equation}\label{Dyn1}
\Delta \mathbf{r}_t := \mathbf{r}_{t+ \Delta t} 
- \mathbf{r}_t 
= \mathbf{ps} \,\Delta w_t  + \mathbf{trd} \, \Delta t, 
\end{equation}
where $ \mathbf{trd} $ is 
a $ k $-dimensional vector which describes the 
deterministic "trend".
Then we have the following (It\^o type) Taylor approximation 
that describes the random evolution of the present value:
\begin{equation}\label{Taylor1}
\begin{split}
\frac{\Delta (\mathrm{PV})}{\mathrm{PV}} 
& :=  \frac{\mathrm{PV} ( \mathbf{r} + \Delta 
\mathbf{r}_t ; t + \Delta t ) -
\mathrm{PV} (\mathbf{r}; t) }
{\mathrm{PV}} \\
&\simeq - D_{\mathbf{ps}} \,\Delta w_t + 
\frac{1}{2} D_{\mathbf{ps}}^2 \, (\Delta w_t )^2 \\
& \qquad \quad + \left\{ \frac{\partial_t (\mathrm{PV})}{\mathrm{PV}} 
-D_\mathbf{trd} \right\}\,\Delta t,
\end{split}
\end{equation}
where $ D_\mathbf{trd} $ is defined in the same way as (\ref{der1}).

From a viewpoint of risk management, 
the expression (\ref{Taylor1}) implies that, 
if we could choose $ \mathbf{CF} $ 
so as to have $ D_{\mathbf{ps}} = 0 $
(and $ D^2_{\mathbf{ps}} = 0 $),
then for the short time period the portfolio is riskless.
Roughly, the duration corresponds to so called {\em delta} and 
the convexity, {\em gamma}. 

These criteria, however, work well only if 
the movements of the term structure of interest rate
are limited to parallel shifts, 
as pointed out in the last part of Chapter 4 of \cite{Hull}, 
which is quoted above. 
To have richer descriptions, 
we want to make the dynamics {\bf multi-factor}:
typically we want to assume
\begin{equation}\label{Dyn2}
\Delta \mathbf{r}_t = 
\sum_{j=1}^n \mathbf{x}^j  \,\Delta w^j_t 
+ \mathbf{trd} \, \Delta t, 
\end{equation}
where $ \mathbf{x}^jj $'s are $ k $-dimensional vectors and 
$ \Delta w^j_t $'s are real valued random variables
with 
\begin{equation}\label{cov}
\Ex [\Delta w^j_t] = 0 \,\,\text{for all $ j $, and}\,\,
\mathrm{Cov} (\Delta w^i_t, \Delta w^j_t ) = 
\begin{cases}
\Delta t & \text{if  $ i=j $,} \\
0 & \text{otherwise}.
\end{cases}
\end{equation}
Here $ n $ is the number of {\em factors}. 
If we write 
$ \Delta \mathbf{w}:= (\Delta w^1,...,\Delta w^n) $, 
the equation (\ref{Dyn2}) is rewritten as
\begin{equation*}
\Delta \mathbf{r}_t = 
\langle \mathbf{x}, \Delta \mathbf{w}_t \rangle 
+ \mathbf{trd} \, \Delta t,
\end{equation*}
or coordinate-wisely
\begin{equation*}
\Delta \mathbf{r}_t(T_l) = 
\langle \mathbf{x}(T_l), \Delta \mathbf{w}_t \rangle 
+ \mathbf{trd}(T_l) \, \Delta t,\,\,l=1,...,k.
\end{equation*}

Define $ D_{\mathbf{x}_l} $'s 
in the same manner as (\ref{der1}), which 
should be called {\em generalized durations}, and define also
\begin{equation}\label{Taylor22}
D^2_{\mathbf{x}^i, \mathbf{x}^j } := (\mathrm{PV})^{-1}
\partial_{\mathbf{x}^i} \partial_{\mathbf{x}^j} (\mathrm{PV})
= (\mathrm{PV})^{-1}
\sum_{l=1}^k \mathbf{x}^i (T_l) \mathbf{x}^j (T_l) (T_l - t)^2 e^{-\mathbf{r}(T_l) (T_l-t)},
\end{equation}
which would be called {\em generalized convexities}.
In these settings (\ref{Taylor1}) turns into
\begin{equation}\label{Taylor2}
\begin{split}
\frac{\Delta (\mathrm{PV})}{\mathrm{PV}} 
& \simeq - \sum_{j=1}^n D_{\mathbf{x}^j} ( \Delta w^j_t )
+ \frac{1}{2} \sum_{ 1 \leq i, j\leq n} D^2_{\mathbf{x}^i, 
\mathbf{x}^j} 
( \Delta w^i_t \Delta w^j_t  ) \\
& \hspace{4cm}+ 
\{ (\mathrm{PV})^{-1} 
\partial_t (\mathrm{PV}) -D_\mathbf{trd} \}\,\Delta t \\
& \text{(in the vector notation with 
$ \mathbf{D}_{\mathbf{x}} = (D_{\mathbf{x}}^1,...,D_{\mathbf{x}}^n) $)} \\
& =  - \langle \mathbf{D}_{\mathbf{x}}, \Delta \mathbf{w}_t \rangle
+ \frac{1}{2PV_t} \langle 
\mathbf{D}_{\mathbf{x}} \otimes \mathbf{D}_{\mathbf{x}},
\Delta \mathbf{w}_t \otimes \Delta \mathbf{w}_t \rangle \\
& \hspace{4cm}+ 
\{ (\mathrm{PV})^{-1} 
\partial_t (\mathrm{PV}) -D_\mathbf{trd}\}\,\Delta t \\
&\text{(in this notation convexities are not needed!)}.
\end{split}
\end{equation}

The parameters $ \mathbf{x}^j $'s 
(including the number of factors $ n $)
can be determined, for example, 
by the {\em principal component analysis} (PCA) 
(see \cite{Li-Sch} or \cite[Section 6.4]{CLM}, 
and in the present paper 
a brief survey of PCA is given in Appendix \ref{PCA}.).
It should be noted, however, 
that such a statistical estimation method
is founded on the basis of {\bf stationarity}:
one needs plenty of homogeneous sample data.

\subsection{Stationary models}\label{forward00}
T.S.Y. Ho, and S. B. Lee \cite{HL} 
proposed such an arbitrage-free {\em stationary}
interest rate model.
In essence Ho-Lee model assumes that 
\begin{equation}\label{paramet}
\{ P^{(t)} (T) 
:= - T^{-1} \log \mathbf{r}_t ( t+T):T \in \Z_+ \} 
\end{equation}
is a {\em random walk}
in $ \R^{\Z_+} $. 
By a {\em random walk}
we mean a sum of i.i.d. random variables.
In other words, a random walk is a 
process with stationary independent 
increments, and hence is consistent with 
PCA.
Actually in Ho-Lee model increments of the random walk 
is coin-tossing, and hence 
is described by a binomial {\em recombining tree}
(see Fig \ref{rct1} and Fig \ref{TStree} in section \ref{revisited}). 

It is widely recognized, however, that 
the Ho-Lee model is insufficient in that it 
can describe the movements of 
``parallel shift" only (see e.g. the quotation in the beginning of
section \ref{revisited}).
In sections \ref{revisited} and \ref{MD1} we will study this 
widely spread belief in a careful way 
and show that this is not because 
the model is one-factor but 
in fact the problem lies 
in its parameterization of (\ref{paramet}) 
(Theorem \ref{multi1}).

In section \ref{MFG}, we will reveal that
this puzzle is resolved by shifting to 
{\em forward rates}\footnote{
A forward rate is a pre-agreed rate 
for borrowing during a pre-agreed future time interval.}. 
To exclude trivial arbitrage opportunities 
continuously compound forward rate over $ (T, T'] $ 
at time $ t $ should be given by
\begin{equation}\label{forward1}
F_t(T,T') := \frac{\mathbf{r}_t (T') (T'-t) 
- \mathbf{r}_t (T) (T-t) }{T' -T}.
\end{equation}
In this shift 
it is assumed that for every $ T $, 
the increments $ \Delta F_t (T) := F_t (T) - F_{t-\Delta t} (T) $
of the instantaneous forward rate 
$ F_t (T) := F_t (T,T+\Delta t) $ are i.i.d..,
or more precisely, we assume that for a sequence of 
i.i.d. random variables $ \Delta \mathbf{w} $ satisfying (\ref{cov}), 
there exists $ \sigma (T) \in \R^n $ and $ \mu (T) \in \R $ 
for each $ T $ 
such that for arbitrary $ 0 \leq t \leq T $
\begin{equation}\label{forward201}
\Delta F_t (T) = \langle \sigma (T),  \Delta \mathbf{w}_t \rangle 
+ \mu (T) \Delta t. 
\end{equation}
The stochastic processes $ F_t (T) $'s are {\bf stationary} in the sense 
$ \sigma (T) $ and $ \mu (T) $ are constant with respect to 
the time parameter $ t $. 

In sharp contrast with the cases of 
Ho-Lee and its direct generalizations 
there always exists an 
arbitrage-free stationary model 
for given term structure of volatilities 
$ \sigma (\cdot) : \Z_+ \Delta t
 \to \R $. 
In fact, we have the following.
\begin{theorem}\label{consistency}
For any given (estimated) term structure of volatilities 
$ \sigma (\cdot) :\Z_+ \Delta t \to \R $
there always exists an arbitrage-free 
{\bf stationary} forward rate model 
given by (\ref{forward201}).
Precisely speaking, the following holds.
The model with (\ref{forward201})
is arbitrage-free if
\begin{equation}\label{driftcond1}
{\mu} (T)
=\frac{1}{\Delta t} 
\log \frac{\Ex[\exp { \langle  \rho (T), \Delta \mathbf{w} \rangle} ] }
{\Ex[\exp { \langle  \rho (T+ \Delta t), \Delta \mathbf{w} \rangle} ]}, 
\end{equation}
where  $ \rho (T) = \Delta t \sum_{ 0 \leq u \leq T} \sigma (u) $. 
\end{theorem}
A proof of Theorem \ref{consistency} will be given in section \ref{MFG}. 

\ 

Though the shift to forward rate model
parallels with the celebrated works by
D. Heath, R.A. Jarrow and A. Morton \cite{HJM1} and \cite{HJM2},
our model presented here is an innovation in that 
\begin{enumerate}[(i)]
\item The model is stationary, and hence consistent with PCA. 
(It is not so clear whether the drift $ {\mu} $
can be constant in time under no-arbitrage hypothesis.) 
\item It is a generalization of Ho-Lee's binomial model 
in the sense that it is within a discrete time-state framework and 
it can be described by 
a {\bf recombining tree}; 
If $ \Delta w^j $'s have only finite possibilities
(in fact we will construct under the assumption of 
$(a2') $ in section \ref{MD1}.)
then the tree (=all the scenarios) of 
$ \mathbf{w} \equiv (w^1,...,w^n) $, 
hence the tree of $ \mathbf{r} $, becomes
{\bf recombining} since 
$ \Delta \mathbf{w} (t) (w) + \Delta \mathbf{w} (t') (w') $
equals $ \Delta \mathbf{w} (t') (w) 
+ \Delta \mathbf{w} (t) (w') $ for any given time
$ t, t' $ and one step scenarios $ w, w' $.

\item It is a {\bf multi-factor} model; along the line 
of (\ref{Taylor2}) and (\ref{Taylor22}) 
it gives an alternative sensitivity analysis beyond {\em duration}
and {\em convexity} (see section \ref{ASA}).  
\item Furthermore, 
this class can be seen as a discrete-time analogue of
(a special cases of)
multi-factor generalization 
of Ritchken-Sankarasubramanian's model \cite{RS}
due to K. Inui and M. Kijima \cite{IK} (see section \ref{PSTL}). 
\end{enumerate}


\section{Linear models}\label{linear}

\subsection{Ho-Lee model revisited}\label{revisited}
\begin{quote}
The Ho-Lee model was a significant improvement over 
what came before, 
but it had a number of failings. 
It was a one-factor model, 
and the way the term structure evolved over time 
was through parallel shifts. (R. Jarrow \cite{Jar})
\end{quote}
Let $ P^{(t)} (T) $
denote the price at time $ t $ of 
zero-coupon bond whose time to maturity is $ T $:
\footnote{Note that this notation 
is different from the standard one. 
The parameter $ T $ usually stands for 
the maturity date as we have done and will do
for $ \mathbf{r} $ and others.
The parameterization
used here originates from \cite{HL}.
}
To avoid confusions, below we give a reversed form of 
(\ref{paramet}): 
\begin{equation}\label{Ptor}
\mathbf{r}_t (T) = \frac{1}{T-t} \log P^{(t)} (T-t), \quad (T>t).
\end{equation}
Here $ t $ and $ T $ are non-negative integers; in particular,  
$ \Delta t $ is set to be $ 1 $. 
Note that $ P $ must be positive at least and 
$ P^{(t)} (0) = 1 $ without fail.

The Ho-Lee model is characterized by the following assumptions 
$(a1)$ to $(a4)$. 

\begin{enumerate}
\item[$(a1)$] The random variable $ P^{(t)} (T) $ depends only on
so-called {\em state variable} $ i  = i_t $ [(A4) of pp1013 in \cite{HL}].
\end{enumerate}
\begin{enumerate}
\item[$(a2)$] For each $ t $, $ \Delta i_t = i_{t+1} - i_t $ 
can take only two values; as in \cite{HL}
the two are $ 1 $ and $ 0 $
[(4) of pp1014 in \cite{HL}]. 
\end{enumerate}
As a consequence of $(a2)$, 
the state space of $ i $ is $ \Z $,
and if we take $ i_0 =0 $, then it can be $ \Z_+ $.
We sometimes write $ P^{(t)}_i (T) $, and in this case
we think of $ P $ as a function on $ \Z_+ \text{(time)} \times 
\Z_+ \text{(state)} \times \Z_+  \text{(maturity)} $.

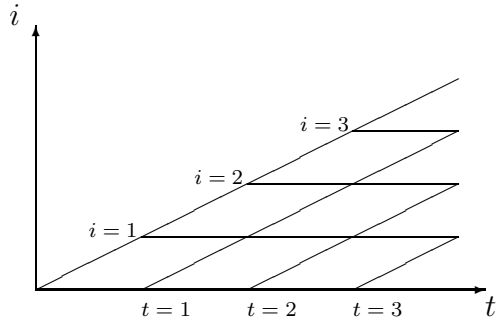
\begin{figure}[h]
\begin{center}
\begin{picture}(170, 100)(-20,0)
\put(0,0){\vector(1,0){170}}
\put(0,0){\vector(0,1){100}}
\put(0,0){\line(2,1){160}}
\put(40,20){\line(1,0){120}}
\put(80,40){\line(1,0){80}}
\put(120,60){\line(1,0){40}}
\put(40,0){\line(2,1){120}}
\put(80,0){\line(2,1){80}}
\put(120,0){\line(2,1){40}}
\put(170,-10){$ t $}
\put(-10,100){$ i $}
\put(40,-10){\scriptsize $ t=1 $}
\put(80,-10){\scriptsize $ t=2 $}
\put(120,-10){\scriptsize $ t=3$}
\put(20,20){\scriptsize $i=1$}
\put(60,40){\scriptsize $ i=2$}
\put(100,60){\scriptsize $ i=3$}
\end{picture}
\caption{Recombining tree \#1}\label{rct1}
\end{center}
\end{figure}

\begin{figure}[h]
\begin{center}
\begin{picture}(200,120)(0,-10)
\put(20, 20){\scriptsize $ P^{(0)}_0 (\cdot) $}
\put(0,0){\vector(1,0){45}}
\put(0,0){\vector(0,1){30}}
\put(0,20){\line(5,-1){10}}
\put(10,18){\line(5,-2){10}}
\put(20,14){\line(5,-1){10}}
\put(30,12){\line(1,0){10}}
\put(50, 15){\vector(1,0){10}}
\put(90, 20){\scriptsize $ P^{(1)}_0 (\cdot) $}
\put(70,0){\vector(1,0){45}}
\put(70,0){\vector(0,1){30}}
\put(70,20){\line(5,-1){10}}
\put(80,18){\line(5,-2){10}}
\put(90,14){\line(5,-1){10}}
\put(100,12){\line(1,0){10}}
\put(50, 40){\vector(2,1){10}}
\put(90, 60){\scriptsize $ P^{(1)}_1 (\cdot) $}
\put(70,40){\vector(1,0){45}}
\put(70,40){\vector(0,1){30}}
\put(70,60){\line(5,-2){10}}
\put(80,56){\line(5,-1){10}}
\put(90,54){\line(5,-1){10}}
\put(100,52){\line(1,0){10}}
\put(120, 15){\vector(1,0){10}}
\put(160, 20){\scriptsize $ P^{(2)}_0 (\cdot) $}
\put(140,0){\vector(1,0){45}}
\put(140,0){\vector(0,1){30}}
\put(140,20){\line(5,-1){10}}
\put(150,18){\line(5,-2){10}}
\put(160,14){\line(5,-1){10}}
\put(170,12){\line(1,0){10}}
\put(120, 60){\vector(1,0){10}}
\put(120, 40){\vector(2,1){10}}
\put(120, 80){\vector(2,1){10}}
\put(160, 60){\scriptsize $ P^{(2)}_1 (\cdot) $}
\put(140,40){\vector(1,0){45}}
\put(140,40){\vector(0,1){30}}
\put(140,60){\line(5,-2){10}}
\put(150,56){\line(5,-1){10}}
\put(160,54){\line(5,-1){10}}
\put(170,52){\line(1,0){10}}
\put(190, 60){\vector(1,0){10}}
\put(190, 70){\vector(2,1){10}}
\put(190, 20){\vector(1,0){10}}
\put(190, 30){\vector(2,1){10}}
\put(140, 90){\scriptsize $ P^{(2)}_2 (\cdot) $}
\put(170, 90){\vector(1,0){10}}
\put(170, 100){\vector(2,1){10}}
\put(190, 100){$\cdots$}
\put(210, 70){$\cdots$}
\put(210, 40){$\cdots$}
\put(210, 20){$\cdots$}
\end{picture}
\end{center}
\caption{Term structures attached to the recombining tree}\label{TStree}
\end{figure}
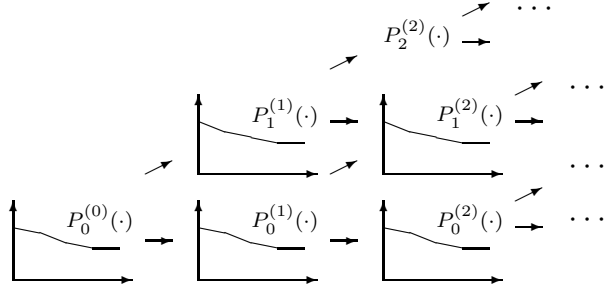

The {\em tree} of the assumption (a2) 
is illustrated in Fig \ref{rct1}.
The first assumption (a1)
also says that the function (often referred to as {\em term structure}) 
$ P^{(t)} ( \cdot ) : \Z_+ \to \R $,
or equivalently the random vector 
$ ( P^{(t)} (1), P^{(t)}(2) $,
$...$, 
$P^{(t)} (T),...) $
depends only on $ i_t $. In other words, 
a realization of term structure is 
attached to each node of the tree of 
Fig \ref{rct1} (See Fig \ref{TStree}). 

The next assumption $(a3)$ makes Ho-Lee model 
as it is [(7), (8), and (10) of pp1017 in \cite{HL}]. 
\begin{enumerate}
\item[$(a3)$] The model is {\em stationary}
in the sense that 
$ P_{i_{t}}^{(t)} (T) P^{(t-1)}_{i_{t-1}} (1) / P_{i_{t-1}}^{(t-1)} (T+1) $
depends only on $ T $ and $ \Delta i_t =i_{t}- i_{t-1} $. 
\end{enumerate}
The assumption $ (a3) $ says that 
if $ \{ \Delta i_t : t \in \N\} $ are i.i.d., 
so are $  \{ P_{i_{t}}^{(t)} (\cdot ) P^{(t-1)}_{i_{t-1}} 
(1) $ $/ P_{i_{t-1}}^{(t-1)} ( \cdot +1) : (t \in \N) \} $. 
In this sense we call it stationary.
Actually we implicitly assume that under the real world measure 
$ \Delta i_t $ are i.i.d..
Denote
\begin{equation}
H (T, \Delta i_t ):= 
\frac{P_{ i_{t}}^{(t)} (T) P^{(t-1)}_{i_{t-1}} (1)}
{ P_{i_{t-1}}^{(t-1)} (T+1) }.
\end{equation}
Recall that 
$ h(T) := H (T, 1) $ and $ h^* (T) := H(T,0) $ are
called perturbation function(s) in \cite{HL}.

The next condition $(a4)$ says that 
the model is arbitrage-free (see Theorem \ref{NAS} in Appendix).
\begin{enumerate}
\item[$(a4)$] There exists $ \pi \in (0,1) $ 
such that
for arbitrary $ t, i \in \Z_+ $ and $ T > 0 $, 
\begin{equation*}
P_{i}^{(t)} (T+1) =  \pi P^{(t)}_i (1) P_{i+1}^{(t+1)} (T) 
+ (1- \pi) P^{(t)}_i (1) P_{i}^{(t+1)} (T). 
\end{equation*}
\end{enumerate}

Under these assumptions, 
dynamics of the entire term structure 
is fully determined.
More precisely, we have the following.
\begin{theorem}\label{holee}
(i) Under the assumption $(a1)$, $(a2)$, and $(a3)$, 
$ \{ P^{(t)} (\cdot): t \geq t_0 \} $ as $ \R^{\Z_+} $ valued 
stochastic process is parameterized by given initial term structure
$ \mathbf{r}_{t_0} $, the constant volatility $ \sigma $,
and the drift function $ \mu : \Z_+ \to \Z_+ $ as 
\begin{equation}\label{solution00}
P^{(t)}_i (T) = \frac{ P^{(0)} (t+T)}{P^{(0)} (t)} 
\exp \{ - T \sigma (2 i - t)
- \sum_{u=0}^{T-1}\sum_{v=1}^{t} \mu (u+v-1) \}.  
\end{equation}
(ii) The no-arbitrage condition $(a4)$ is equivalent to 
the following ``drift condition": for a $ \pi \in (0,1) $ we have
\begin{equation}\label{drift1}
\mu (T) = 
\log \frac{ \pi  e^{-(T+1) \sigma } + (1-\pi) e^{(T+1) \sigma} }
{\pi e^{ - T \sigma } + (1-\pi) e^{T \sigma} }.
\end{equation} 
\end{theorem} 
This theorem will be proved as a special case 
($ n =2 $ and $ i = 2 w -1 $) of 
Theorem \ref{multi1} below. 
\begin{coro}[ Ho and Lee \cite{HL}, (19) and (20) of pp 1019 ] 
Under $ (a1) $, $ (a2) $, and $ (a3) $, 
no-arbitrage condition $ (a4) $ is equivalent to 
the existence of $ \pi \in (0,1) $ and $ \delta > 0 $ such that 
\begin{equation}\label{perturb}
h (T) = \frac{1 }
{\pi  +  (1- \pi) \delta^T} \,\,\text{and}\,\,
h^* (T) = \frac{ \delta^T }
{\pi  +  (1- \pi) \delta^T}
\end{equation}
for all $ T \in \Z_+ $. 
\end{coro}
\begin{proof}
By  (\ref{mutoH}) below we have 
\begin{equation*}
\begin{split}
H (T, \Delta i) & = e^{-  T \sigma  (2\Delta i -1) } 
e^{-  \sum_{u=0}^{T-1} \mu(u) }, \\
& (\text{by substituting (\ref{drift1})} ) \\
& =  e^{- \sigma T (2\Delta i -1) } 
\prod_{u=0}^{T-1} \frac{ \pi  e^{-u\sigma } + (1-\pi) e^{u \sigma}}
{ \pi e^{-(u+1) \sigma } +  (1- \pi)  e^{(u+1) \sigma} } \\
&= \frac{e^{- \sigma T (2\Delta i -1) } }
{\pi e^{-T \sigma } +  (1- \pi)  e^{T \sigma} }.
\end{split}
\end{equation*}
Therefore, we have 
\begin{equation*}
h (T) = H( T,1) =  \frac{e^{- \sigma T  } }
{\pi e^{-T \sigma } +  (1- \pi)  e^{T \sigma} }
= \frac{1 }
{\pi  +  (1- \pi)  e^{2 T \sigma} }
\end{equation*}
and 
\begin{equation*}
 h^* (T) = H(T,0) =  \frac{e^{ \sigma T  } }
{\pi e^{-T \sigma } +  (1- \pi)  e^{T \sigma} }
=  \frac{e^{ 2 \sigma T  } }
{\pi +  (1- \pi)  e^{ 2T \sigma} }.
\end{equation*}
By putting $ \delta = e^{2\sigma} $ we obtain (\ref{perturb}).
\end{proof}

\subsection{Multi-dimensional generalization}\label{MD1}

As we have discussed in section \ref{sensitivity},
by a multi-factor generalization
we want to mean such as (\ref{Dyn2})
and (\ref{forward201}), 
with the moment condition (\ref{cov}).

For the time being let us consider the condition (\ref{cov}).
Suppose that $ \Delta \mathbf{w} $ is defined on a finite set, say, $ S
= \{ 0,1,...,s\} $.
Let $ L (S) := \{ f : S \to \R \} $. We endow the scholar product 
with $ L (S) $ by $ \langle x, y \rangle = \Ex [xy] $.   
Since $ \Delta t =1 $ here in this section,
the condition (\ref{cov}) is rephrased as:
\begin{equation*}
\text{$ \{ 1, \Delta w^1,...,\Delta w^n\} $ is an orthonormal system 
of $ L(S) $}.
\end{equation*}
(For general $ \Delta t $, it is replaced with 
$ \{ 1, \Delta w^1 /\sqrt{\Delta t},...,
\Delta w^n /\sqrt{\Delta t} \} $.)
From this observation
we notice that $ s +1 = \sharp S $ should be greater than $ n $.
If $ s > n $, then we can extend 
$ \{ 1, \Delta w^1,...,\Delta w^n\} $ 
to have an orthonormal basis of $ L (S) $
by constructing (dummy) random variables
$ \Delta w^{n+1},...,\Delta w^{s } $. 
Thus for simplicity we want to assume $ n= s $.

Let us go a little further. Let 
\begin{equation}\label{constrcution1}
m_{i,j} = \Delta w^i (j) \sqrt{\Pr (\{s_j\}) }, \quad 0 \leq i \leq n, 
\, 0 \leq j \leq s,
\end{equation}
where $ \Delta w^0 \equiv 1 $. 
When $ n= s $, then 
$ M = (m_{i,j})_{0 \leq i,j \leq n=s} $
is an orthogonal matrix; following a standard notation,
we have $ M \in O (n+1) $.
Conversely, given an orthogonal matrix 
$ M = (m_{i,j})_{\underline{0} \leq i,j \leq n} \in O (n+1) $
with $ m_{0, j} > 0 $ for all $ j $, 
we can construct a random variable $ \Delta \mathbf{w}
= (\Delta w^1,...,\Delta w^n ) $ satisfying (\ref{cov}) by
setting $ m_{0, j} = \sqrt{\Pr (\{s_j\}) } $ and 
the relation (\ref{constrcution1}). 

With the above considerations in mind, 
we generalize the assumption $ (a2) $
as follows.
\begin{enumerate}
\item[$(a2')$] The state variable $ i_t $ is $ \R^{n+1} $-valued,
and $ \Delta i_t = ( \Delta t, \Delta w^1_t,...,\Delta w^n_t ) $,
where $ ( \Delta w^1_t$,
$...,\Delta w^n_t) $, $ t=1,2,... $ are 
independent
copies of 
$\Delta \mathbf{w} $;
a random variable 
satisfying (\ref{cov}) and taking only $ n+1 $ distinct values. 
\end{enumerate}
Note that $ (a2) + (a1) $ can be seen as a special case of 
$ (a2') + (a1) $ by setting 
$ n=1 $, 
$ \Delta w^1 (0) =-1 , \Delta w^1 (1) = 1 $, and 
$ \Pr (\{0\}) = \Pr (\{1\}) = 1/2 $
In fact, $ \hat{i}_t = w^1_t/2 + t/2 $, 
is a function of 
$ i_t = ( t, w_t ) $, 
and $ \Delta \hat{i}_t $ can be either $ 1 $ or $ 0 $. 
Conversely, starting from $ (a2) $,
the function $ (1,2i -1) $ on $ S $ can be 
orthonormal basis of $ L (S) $. 

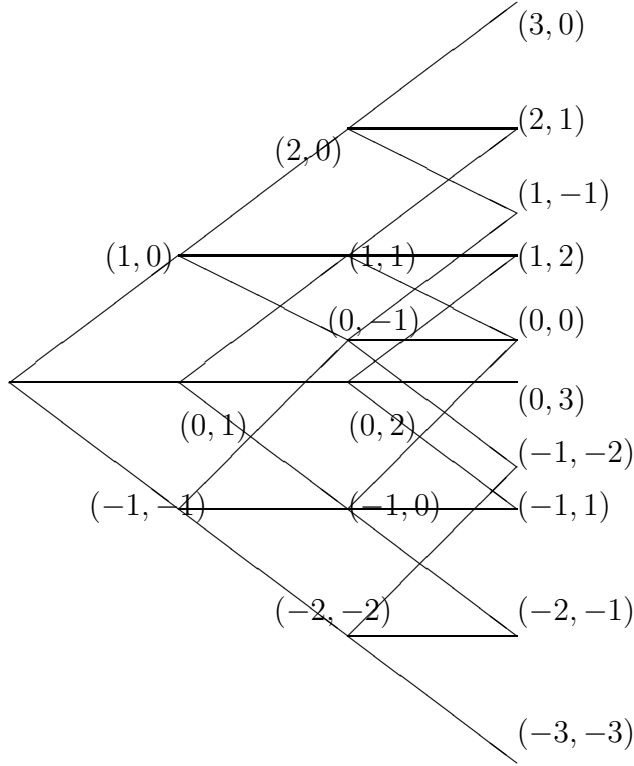
\begin{figure}[h]
\begin{center}
\begin{picture}(340, 280)(-50,-140)
\put(0,0){\line(1,0){64}}
\put(64,-20){$ (0,1) $}
\put(0,0){\line(4,3){64}}
\put(36,44){$(1,0)$}
\put(0,0){\line(4,-3){64}}
\put(30,-50){$(-1,-1)$}
\put(64,0){\line(1,0){64}}
\put(128,-20){$ (0,2) $}
\put(64,0){\line(4,3){64}}
\put(128,44){$(1,1)$}
\put(64,0){\line(4,-3){64}}
\put(128,-50){$(-1,0)$}
\put(64,48){\line(1,0){64}}
\put(64,48){\line(4,3){64}}
\put(100,84){$(2,0)$}
\put(64,48){\line(2,-1){64}}
\put(120,20){$(0,-1)$}
\put(64,-48){\line(1,0){64}}
\put(64,-48){\line(1,1){64}}
\put(64,-48){\line(4,-3){64}}
\put(100,-90){$(-2,-2)$}
\put(128,0){\line(1,0){64}}
\put(192,-10){$ (0,3) $}
\put(128,0){\line(4,3){64}}
\put(192,44){$(1,2)$}
\put(128,0){\line(4,-3){64}}
\put(192,-50){$(-1,1)$}
\put(128,96){\line(1,0){64}}
\put(192,96){$(2,1)$}
\put(128,96){\line(4,3){64}}
\put(192,132){$(3,0)$}
\put(128,96){\line(2,-1){64}}
\put(192,68){$(1,-1)$}
\put(128,48){\line(1,0){64}}
\put(128,48){\line(4,3){64}}
\put(128,48){\line(2,-1){64}}
\put(192,20){$(0,0)$}
\put(128,-48){\line(1,0){64}}
\put(128,-48){\line(1,1){64}}
\put(128,-48){\line(4,-3){64}}
\put(192,-90){$(-2,-1)$}
\put(128,-96){\line(1,0){64}}
\put(128,-96){\line(1,1){64}}
\put(192,-30){$(-1,-2)$}
\put(128,-96){\line(4,-3){64}}
\put(192,-136){$(-3,-3)$}
\put(128,16){\line(4,3){64}}
\put(128,16){\line(1,0){64}}
\put(128,16){\line(4,-3){64}}
\end{picture}
\caption{Recombining tree \#2. 
(Note that the pairs of numbers does not 
mean the values of $ w $.)
}\label{rct2}
\end{center}
\end{figure}

We denote the state space of $ \mathbf{w}_t $ by $\mathfrak{Z} $.
In general the state space of $ i_t = (t, \mathbf{w} ) $
is described by a recombining tree, the number of 
whose nodes grow in an order
of $ t^n $. Detailed studies of the construction of
$ \mathbf{w} $ and its tree are presented in \cite{A3}.
Here we illustrate a tree 
for $ n =2 $ in Fig \ref{rct2}. 

According as this generalization, 
$(a4) $ also needs to be modified.
\begin{enumerate}
\item[$(a4')$] There exists 
a probability (measure) $ \pi_i (\cdot) \in (0,1) $  on $ S $ 
such that
for arbitrary $ T > 0 $ and $ i = (t,w) \in \mathfrak{Z} $,
\begin{equation*}
P_{i}^{(t)} (T+1) = P^{(t)}_i (1) 
\sum \pi_i (s) P_{i + \Delta i (s)}^{(t+1)} (T). 
\end{equation*}
\end{enumerate}

As we have remarked in section \ref{forward00}, 
the stationarity assumption $(a3)$
is yet too restrictive even under these relaxed conditions.

\begin{theorem}\label{multi1}
(i) Under the assumption $(a1)$, $(a2')$, and $(a3)$, 
dynamics of the entire term structure
$ \{ P^{(t)} (\cdot): t \geq 0 \} $ as $ \R^{\Z_+} $ valued 
stochastic process is parameterized by given initial term structure
$ P^{(0)} (\cdot) :\R^{\Z_+} $, 
the constant volatility $ \sigma  \in \R^n $,
and the drift function $ \mu : \Z_+ \to \R $ as 
\begin{equation}\label{solution0}
P^{(t)} (T) = \frac{ P^{(0)} (t+T)}{P^{(0)} (t)} 
\exp \{ - T \langle \sigma , \mathbf{w}_t - \mathbf{w}_0 
\rangle - \sum_{u=0}^{T-1}\sum_{v=1}^{t} \mu (u+v-1) \}.  
\end{equation}
(ii) The no-arbitrage condition $(a4')$ is equivalent to 
the following ``drift condition": for a probability $ \pi $ on $ S $ 
we have  
\begin{equation}\label{drift2}
\mu (T) = 
\log \frac{ \sum_s 
\pi (s) e^{-(T+1) \langle \sigma , \Delta \mathbf{w} (s) \rangle } }
{\sum_s 
\pi (s) e^{-T \langle \sigma , \Delta \mathbf{w} (s) \rangle }}.
\end{equation}
\end{theorem} 
\begin{proof} 
(i) 
First, notice that the assumption $(a1)$
in terms of $ \mathbf{f}_t (T) 
= F_t ( t+T, t+ T+\Delta t) $, 
the forward rate over $ (T-t, T+1-t] $ at time $ t $,
is stable: 
If $(a1)$ holds, then $ \mathbf{f}_t (T) $ depends only on 
$ i_t $ and vise versa. 
To see this, just use the following expressions:
\begin{equation*}
\mathbf{f}_t (T) =  \log\frac{P^{(t) }(T)}{P^{(t)}(T+1)},
\end{equation*}
and 
\begin{equation}\label{converse1}
P^{(t) }(T) = \exp \{ - \sum_{u=0}^{T-1} \mathbf{f}_t (u) \}.  
\end{equation}

Next, observe that the assumption $(a3)$  
in terms of $ \mathbf{f}_t (T) $
turns into the following equation:
\begin{equation}\label{BGMdisc}
\mathbf{f}_{t} (T) - \mathbf{f}_{t-1} (T+1)
= \log H (T, \Delta i_{t}) - \log H (T+1, \Delta i_{t}).
\end{equation}
Thanks to the 
``completeness" assumption, 
the right-hand-side of (\ref{BGMdisc})
is reparametrized as  
$ \langle \sigma (T),  \Delta {i}_{t} \rangle_{\R^n} $
by some $ \sigma (T) \in \R^{n+1} $ for each $ T \in \Z_+ $.
Note that the choice of $ \sigma (T) $ is unique. 
With this representaion, 
the solution to (\ref{BGMdisc}) is obtained as
\begin{equation}\label{solution1}
\mathbf{f}_{t} (T) = \mathbf{f}_0 (T+t) + \sum_{u=1}^{t}
\langle \sigma (T+t-u) , \Delta i_{u} \rangle.   
\end{equation}

Since $ i_t $ is an exchangeable in the sense that 
\begin{equation*}
i (s_1,...,s_t) 
=  i (s_{\tau(1)},...,s_{\tau(t)})
\end{equation*}
for arbitrary permutation $ \tau \in \mathfrak{S}_t $, 
so is $ \mathbf{f}_{t} (T) \equiv \mathbf{f}_t (T; s_1,...,s_t) $,
$ (s_1,...,s_t) \in S^t  $.
Therefore, from the explicit expression (\ref{solution1}) we have
\begin{equation*}
\begin{split}
\sum_{u=1}^{t}
\langle \sigma (T+t-u) , \Delta i (s_u)\rangle
& = \sum_{u=1}^{t}
\langle \sigma (T+t-u) , \Delta i ( s_{\tau(u)} )\rangle \\
&= \sum_{u=1}^{t}
\langle \sigma (T+t-\tau^{-1} (u)) , \Delta i ( s_u)\rangle
\end{split}
\end{equation*}
for arbitrary permutation $ \tau $. 
In particular, the relation holds for any transposition.
Therefore, we have
\begin{equation*}
\begin{split}
& \langle \sigma (T+t-u) -\sigma (T+t-u'), \Delta i (s_u)
- \Delta i (s_{u'}) \rangle_{\R^{n+1}} \\
& \quad 
\text{(since the first coordinate of $ i $ is independent of $ s $)} \\
& \quad = 
\langle \sigma' (T+t-u) -\sigma' (T+t-u'), 
\Delta \mathbf{w} (s_u)
- \Delta \mathbf{w} (s_{u'}) \rangle_{\R^n},
\end{split}
\end{equation*}
where $ \sigma' $ is the $ n $-dimensional vector 
projected from the latter $ n $ components of $ \sigma $. 

On the other hand, since $ M = (m_{i,j})_{ 0 \leq i,j \leq n} $ 
is regular, so is (at least) one of 
the $ n \times n $ matrices $ M_k :=
(m_{i,j})_{ 1 \leq i \leq n, 0 \leq j \leq n, j \ne k  }$, $ k=0,1,...,n $.
This is because $ \det M = \sum_k m_{0,k} (-1)^{n+k} \det M_k $.
Consequently, setting $ s_{u'} = k $, 
we have that $ \{ \Delta \mathbf{w} (j)
- \Delta \mathbf{w} (k) $, $ j=0,1,...,n $, 
$ j \ne k \} $ 
span the whole space $ \R^n $.
Therefore, we have 
\begin{equation}\label{flatter}
\text{$ \sigma' (T+t-u) -\sigma' (T+t-k) = 0 $
for arbitrary $ 1 \leq u  \leq t $. 
}
\end{equation}
Since (\ref{flatter}) is true for any $ T \in \Z_+ $, 
$ \sigma' (T) $ equals to
a constant $ \sigma  \in \R^n $ irrespective of $ T $. 

Reminding that first component of $ \sigma $ 
can be taken arbitrary, which we denote by $ \mu $,
we have 
\begin{equation*}
\mathbf{f}_{t} (T) = \mathbf{f}_0 (T+t) + 
\langle \sigma , \mathbf{w}_t - \mathbf{w}_0 
\rangle + \sum_{u=1}^{t} \mu (T+u-1).
\end{equation*}
Hence, by (\ref{converse1}) 
we have (\ref{solution0}).

(ii) Back to the original parameterization, we have seen 
\begin{equation*}
\log H (T, \Delta i_{t}) - \log H (T+1, \Delta i_{t}) 
= \langle \sigma , \Delta \mathbf{w}_t \rangle 
+ \mu (T), \quad T \in \Z_+.
\end{equation*}
Therefore, (since $ H (0, \cdot ) = P (0) = 1 $)
\begin{equation}\label{mutoH}
\log H (T, \Delta i ) = -T \langle \sigma , \Delta \mathbf{w}_t \rangle 
- \sum_{u=0}^{T-1} \mu (u). 
\end{equation}
On the other hand, since the assumption $ (a4') $ 
in terms of $ H $ is 
\begin{equation*}
1 = 
\sum_s \pi_i (s) H (T, \Delta i (s) ), 
\end{equation*}
we have
\begin{equation*}
e^{\sum_{u=1}^{T-1} \mu (u)} = \sum_s 
\pi_i (s) e^{-T \langle \sigma , \Delta \mathbf{w}_t (s) \rangle }, 
\end{equation*}
and hence (\ref{drift2}). 
\end{proof}

\begin{remark}
{\em Dynamics of term structure in terms of spot rates is now 
\begin{equation*}
\Delta \mathbf{r} ( \cdot ) 
= \langle \sigma , \Delta \mathbf{w} \rangle + (\mbox{deterministic}).
\end{equation*}
This means that under ($a3$) Ho-Lee model 
is still poor; it can describe only parallel shifts.
In this sense this is {\em not} a multi-factor generalization.
}
\end{remark}

\begin{remark}
{\em (This was suggested by Damir Filipovic in a private communication) 
The equation (\ref{BGMdisc}), or reparametrized one
\begin{equation}\label{BGMd2}
\mathbf{f}_{t} (T) - \mathbf{f}_{t-\Delta t} (T+\Delta t )
= \langle \sigma , \Delta \mathbf{w}_t \rangle 
+ \mu (T), 
\end{equation}
is a discrete-time analogue of
a class of Stochastic Partial Differential Equations (SPDE) 
studied by Brace, Gatarek, and Musiala \cite{BGM}.
In fact,  at least symbolically (\ref{BGMd2}) can be rearranged into
\begin{equation*}
\Delta \mathbf{f}_t (T)  
=  ( \partial_T \mathbf{f}_{t-\Delta t} + \mu' (T) ) \Delta t + 
\langle \sigma , \Delta \mathbf{w}_t \rangle,
\end{equation*}
where $ \mu' (T) = \mu(T) /\Delta t $, 
\begin{equation*}
\Delta \mathbf{f}_t (T)  \equiv 
\mathbf{f}_{t} (T) - \mathbf{f}_{t-\Delta t} (T)
\end{equation*}
and 
\begin{equation*}
\partial_T \mathbf{f}_{t-\Delta t} \equiv
\frac{1}{\Delta t} \{ 
\mathbf{f}_{t-\Delta t} (T+\Delta t ) - \mathbf{f}_{t-\Delta t} (T) \}.
\end{equation*}
}
\end{remark}

\subsection{Multi-{factor} 
generalization}\label{MFG}

\begin{quote}
Instead of focusing upon bond prices as in Ho
and Lee (1986), we concentrate on forward rates. 
This "modification" makes
the model easier to understand and to perform mathematical analysis.
(D. Heath, R. Jarrow and A. Morton \cite{HJM1})
\end{quote}
We have shown in the previous section
that $ (a3) $ is too restrictive to 
have rich models. 
One may well ask if 
there would be another stationary assumption other than $ (a3) $.
The simplest alternative candidate may be 
\begin{equation}\label{simplest}
-\frac{1}{T} \log P^{(t)} (T) =  \langle  \sigma (T), 
\mathbf{w}_t \rangle + \mu(T) \,t ,
\end{equation}
but:
\begin{proposition}
The model (\ref{simplest}) 
cannot be arbitrage-free unless 
$ \sigma (T) = \sigma (1) $ for arbitrary $ T \geq 1 $.
Namely only parallel shifts are allowed for.
\end{proposition}
\begin{proof}
By Theorem \ref{NAS} (iv) we know that 
a model is arbitrage free if and only if 
\begin{equation*}
\Ex^\pi [ P^{(t+\Delta t)}(T) P^{(t)} (1) | \F_t ] = P^{(t)} (T+1) 
\end{equation*}
with a probability $ \pi $ on $ \F_t $. 
Applying this to (\ref{simplest}), we have
\begin{equation*}
\begin{split}
& \Ex^\pi [e^{\langle  \sigma (T), 
\Delta \mathbf{w} \rangle + \mu(T) \Delta t}] \\
& = 
\exp \{ \langle  (T+1)\sigma (T+1)- T\sigma (T) - \sigma (1) , 
\mathbf{w}_t  \rangle \\
& \qquad \qquad \qquad + (T+1) \mu(T+1) -T\mu(T) - \mu(1) \}. 
\end{split}
\end{equation*}
While the left-hand-side is a constant, the right-hand-side is not
unless $$ (T+1)\sigma (T+1)- T\sigma (T) - \sigma (1) = 0. $$
This means $ \sigma (T) \equiv \sigma (1) $ for arbitary $ T $. 
\end{proof}

To summarize, it would be safe to say that
stationarity requirements are difficult to be compatible with 
no-arbitrage restrictions.
As we have mentioned in section \ref{intro1}, however,
we can resolve the problem by 
changing the parameterization of maturity;
assuming (\ref{forward201}) instead of $ (a3) $. 

Let us introduce (or come back to) the following (standard) notation:
\begin{equation*}
P_t^T = P^{(t)} (t+T) \quad t \leq T.
\end{equation*}
(in this case $ T $ means ``maturity", not ``time to maturity").
Then 
\begin{equation}\label{PtoF2}
F_t (T) = \frac{1}{\Delta t} \log \frac{P^T_t}{P^{T+\Delta t}_t}. 
\end{equation}
and
\begin{equation}\label{FtoP2}
P^T_t = \exp\{ -\Delta t \sum_{u=t}^{T-\Delta t} F_t (u) \}.
\end{equation}
Here and hereafter the unit $ \Delta t $ is not necessarily $ 1 $.

By an {\em interest rate model} we mean 
a family of strictly positive stochastic process $\{ P^T_t \}$,
which are, at each time $ t $, a function of
$ (\Delta \mathbf{w}_1 ,..., \Delta \mathbf{w}_t )$.
In a more sophisticated way of speaking, 
they are
adapted to the natural filtration 
of the random walk $ \mathbf{w} $.
We will denote 
$ \F_t:= \sigma ( \Delta \mathbf{w}_1,...,\Delta \mathbf{w}_t ) $. 
Then $ \{ \F_t \} $ is the natural filtration. 

To prove Theorem \ref{consistency}
we rely on the following theorem, 
which is actually a corollary of 
Theorem \ref{NAS}.
\begin{theorem}\label{SPD}
An interest rate model $ \{ P^T_t\} $ is arbitrage-free 
if and only if there exists a strictly positive $ \{\F_t \}$-adapted 
stochastic process 
$ \mathsf{D}_t $ such that 
\begin{equation}\label{SPDrep}
P^T_t = \Ex [\mathsf{D}_T \,|\, \F_t ]/ \mathsf{D}_t 
\end{equation}
holds for arbitrary $ t \leq T $.
\end{theorem}
A proof of Theorem \ref{SPD} will be given in Appendix \ref{NAR}.
As a corollary we have the following {\em consistency condition
with respect to the initial forward rate curve}
for interest rate modeling.
\begin{coro}
An interest rate model $ \{ P^T_t\} $ is arbitrage-free 
if and only if there exists a strictly positive $ \{\F_t \}$-adapted 
stochastic process 
$ \widehat{\mathsf{D}}_t $ with 
$ \Ex[ \widehat{\mathsf{D}}_t ]= 1 $ for all $ t \geq 0 $ 
such that (i) 
\begin{equation*}
P^T_t = \frac{P^T_0 \Ex [ \widehat{\mathsf{D}}_T \,|\, \F_t ]}
{P^t_0 \widehat{\mathsf{D}}_t} 
\end{equation*}
holds for arbitrary $ t \leq T $, or equivalently
(ii) (in terms of the instantaneous forward rate)
\begin{equation*}
F_t (T) = F_0 (T) 
+ \frac{1}{\Delta t} \log \frac{\Ex [ \widehat{\mathsf{D}}_T \,|\, \F_t ]}{
\Ex[\widehat{\mathsf{D}}_{T+\Delta t} \,|\, \F_t ]}
\end{equation*}
holds for arbitrary $ t \leq T $.
\end{coro}
\begin{proof}
Observing $ P^T_0 = \Ex [\mathsf{D}_T ]/ \mathsf{D}_0 $ by setting 
$ t=0 $ in (\ref{SPDrep}), 
the first statement (i) 
can be implied from (\ref{SPDrep}) by putting $ \widehat{\mathsf{D}}_t 
= \mathsf{D}_t / \Ex [ \mathsf{D}_t] $. 
The converse, from (i) to (\ref{SPDrep}), is trivial. 
The equivalence of (i) and (ii) is obtained from 
(\ref{FtoP2}) and (\ref{PtoF2}). 
\end{proof}

Now we are in a position to prove Theorem \ref{consistency}.

\begin{proof}[Proof of Theorem \ref{consistency}]
Let $ \mathsf{D}_t 
:= \exp { \langle \rho (t), \mathbf{w}_t \rangle} $.
Then 
\begin{equation*}
\begin{split}
 \Ex [\mathsf{D}_T | \F_t ]
&= \Ex [ \exp { \langle \rho (T), \mathbf{w}_T 
\rangle } | \F_t] \\
&= \exp { \langle \rho (T), \mathbf{w}_t 
\rangle} \Ex [ \exp { \langle \rho (T), \mathbf{w}_T 
\rangle}] \\
&= \exp { \langle \rho (T), \mathbf{w}_t 
\rangle} \prod_{u=t+\Delta t}^T 
\Ex [ \exp { \langle \rho (T), \Delta \mathbf{w}_u \rangle}] \\
&= \exp { \langle \rho (T), \mathbf{w}_t \rangle}
( \Ex [ \exp { \langle \rho (T), \Delta \mathbf{w} \rangle}])^{(T-t)/\Delta t}.
\end{split}
\end{equation*}
In particular we have 
\begin{equation*}
 \Ex [\mathsf{D}_t ] = \exp { \langle \rho (t), \mathbf{w}_0 \rangle}
( \Ex [ \exp { \langle \rho (t), \Delta \mathbf{w} \rangle}])^{t/\Delta t}.
\end{equation*}
Let 
\begin{equation*}
\widehat{\mathsf{D}}_t := \mathsf{D}_t /\Ex [ \mathsf{D}_t] =
\exp { \langle \rho (t), \mathbf{w}_t - \mathbf{w}_0 \rangle}
( \Ex [ \exp { \langle \rho (T), \Delta \mathbf{w} \rangle}])^{-t/\Delta t}.
\end{equation*}
Then by Theorem \ref{SPD} the model with
\begin{equation}\label{forwardint}
\begin{split}
F_t (T) &:= F_0 (T) 
+  \frac{1}{\Delta t} \log \frac{\Ex [ \widehat{\mathsf{D}}_T \,|\, \F_t ]}{
\Ex[\widehat{\mathsf{D}}_{T+\Delta t} \,|\, \F_t ]} \\
&= F_0 (T) 
+ \frac{1}{\Delta t} 
\langle \rho (T) -\rho (T+\Delta t) , \mathbf{w}_t - \mathbf{w}_0 \rangle \\
& \quad + \frac{T-t}{(\Delta t)^2} 
\log \Ex[\exp { \langle \rho (T), \Delta \mathbf{w} \rangle} ] \\
& \qquad - \frac{T+\Delta t -t}{(\Delta t)^2} 
\log \Ex[\exp { \langle \rho (T+ \Delta t), 
\Delta \mathbf{w} \rangle} ],
\end{split}
\end{equation}
is arbitrage-free. 
Noting that $ \rho (T) -\rho (T+\Delta t) = \sigma (T) \Delta t$,
from (\ref{forwardint}) we can deduce
\begin{equation*}
\Delta F_t (T)
= \langle \sigma (T) ,  \Delta \mathbf{w}_t \rangle 
+ \frac{1}{\Delta t}
\log \frac{\Ex[\exp { \langle \rho (T), \Delta \mathbf{w} \rangle} ] }
{\Ex[\exp { \langle \rho (T+ \Delta t), \Delta \mathbf{w} \rangle} ]}.
\end{equation*}
Thus we have proved the assertion.
\end{proof}

\section{Application to sensitivity analysis}\label{ASA}
\subsection{parameter estimation}
Let us come back to the sensitivity analysis. 
To fit reality, we assume that we can observe only
such forward rates as $ F_t (T_i, T_{i+1})$ (for past $ t $, of course)
for a coarse set of maturities $ 
\{T_1,....,T_k\} $; 
it may happen  that $ T_i - T_{i+1} \gg \Delta t $. 

It should be noted that from (\ref{forward1})
we have 
\begin{equation}\label{forward12}
(T' -T) F_t (T , T' ) + 
(T'' - T') F_t (T', T'') = (T'' - T) F_t (T, T'') 
\end{equation}
for arbitrary $ t <T < T' < T'' $.
Therefore,
the assumption (\ref{forward201}) implies
the following:
\begin{equation}\label{forward401}
\begin{split}
\Delta F_t(T_i,T_{i+1}) & :=  F_{t+\Delta t} (T_i, T_{i+1})
- F_t( T_i, T_{i+1} )  \\
&= \langle {\sigma} (T_i, T_{i+1}) , \Delta \mathbf{w}_t \rangle  
+ {\mu} (T_i, T_{i+1}) \Delta t, 
\end{split}
\end{equation} 
where 
\begin{equation}\label{volatility100}
 {\sigma}_j (T_i, T_{i+1}) 
=  \frac{\Delta t }{T_{i+1}- T_i }
\sum_{u=T_i}^{T_{i+1} - \Delta t} \sigma_j (u)
= \frac{\rho(T_{i+1}) - \rho (T_{i})}{T_{i+1}-T_i},
\end{equation}
for $ j=0,1,...,n $, with the convention that 
$ \sigma_0 = \mu $. 
The drift condition (\ref{driftcond1}) implies 
\begin{equation}\label{driftcond2}
\mu (T_i, T_{i+1}) = \frac{1}{T_{i+1} -T_i} 
\log \frac{\Ex[\exp { \langle  \rho (T_i), \Delta \mathbf{w} \rangle} ] }
{\Ex[\exp { \langle  \rho (T_{i+1}), \Delta \mathbf{w} \rangle} ]}.
\end{equation}
Note that for every $i $, $ \Delta F_t(T_i,T_{i+1} ) $ 
are again i.i.d. for all $ 0\leq t \leq T_i $.

In practice, up to time $ t $ ($< T_1 $) 
we may obtain sufficiently many (though it is less that $ t/\Delta t $)
homogeneous sample data of 
$ \Delta F (T_i, T_{i+1}) $ for $ i = 1,...,k-1 $, 
and then 
by PCA (see \ref{PCA}), 
we can estimate (together with $ n $ itself)
the {\em volatility matrix} $ 
(\sigma_{i,j}) $
with $ \sigma_{i,j} = \sigma_j (T_i, T_{i+1}) $ 
$ \in \R^k \otimes \R^n $.

To obtain 
a full term structure of volatilities $ \sigma : \Z_+ \Delta t \to \R $
one needs to interpolate the volatility 
matrix $ \sigma $. 
An easiest interpolation is given by
\begin{equation}\label{interpolation1}
\sigma^j (T) = \sigma_{ij} \equiv {\sigma}^j (T_i, T_{i+1} )  \,\,
\text{if $ T_i < T\leq T_{i+1}, i=0,1,...,k $}, 
\end{equation}
where $ j=0,1,...,n $
with conventions of $ T_0 = -\infty $, $ T_k = +\infty $,
and $ \sigma_{i 0} = \sigma^0 (T_i, T_{i+1}) = \mu (T_i, T_{i+1}) $.

\begin{remark}{\em
It may happen that a statistical estimator of $ \mu $ 
(possibly by sample average)
does not necessarily satisfy the drift condition.
Note, however, that neither (\ref{driftcond1}) nor (\ref{driftcond2}) 
is if and only if condition; there might be other arbitrage-free models
that are consistent with estimated volatility matrix. 
A careful study on the existence of arbitrage-free model 
that is consistent with both estimated volatility and \underline{average}
\footnote{This is not the one by drift conditions, of course. }
is still open and would be interesting, both theoretically and practically. 
}
\end{remark}

\begin{remark}{\em 
It should be noted that in our model, 
we first work on the real world measure,
and then estimate the risk neutral measure 
(see Theorem \ref{NAS} (iv)).
The stationarity cannot expect in general under
the risk neutral measure. }
\end{remark}

\subsection{A new framework of sensitivity analysis}
Based on the stationary model by Theorem \ref{consistency}
and the estimated parameters we now introduce 
a new framework of sensitivity analysis.

\begin{theorem}
We have the following random expansion with respect to
$ \Delta t $ of the
present value $ \mathrm{PV} $ of a given
cash flow $ \mathbf{CF} $: 
\begin{equation}\label{Taylor8}
\begin{split}
\frac{\Delta (\mathrm{PV})}{\mathrm{PV}} 
& \simeq \sum_l \mathbf{CF} (T_l) e^{-(T_l-t) \mathbf{r}_t (T_l) }
\bigg\{ 
\left\langle \{\rho (T_l) - \rho(t) \}/\mathrm{PV}_t , 
\Delta \mathbf{w}_t \right\rangle \\
& + \frac{1}{2\mathrm{PV}_t} \langle 
\{ \rho (T_l) - \rho(t) \} \otimes \{ \rho (T_l) - \rho(t) \},  
\Delta \mathbf{w}_t \otimes \Delta \mathbf{w}_t \rangle \\
&+\bigg(    
\langle \sigma (T_1, T_2), \mathbf{w}_{t} 
- \mathbf{w}_0 \rangle - \{ \rho^0 (T_l ) -  \rho^0 (t) \}
+ \rho^0(t) t \bigg) \Delta t \bigg\}
+ \mathrm{o} (\Delta t).
\end{split}
\end{equation}
Here by $ \rho^0 (t) $ 
we mean $ \sum_{0 < u \leq t} \mu(u) $.
Consequently, 
\begin{equation}\label{durations-g}
\mathrm{PV}^{-1} \sum_l (\rho^j (T_l)- \rho^j (t))
\mathbf{CF} (T_l) e^{-(T_l-t) \mathbf{r}_t (T_l) },
\,\,j=1,...,n
\end{equation}
are generalized durations and
\begin{equation}\label{convexities-g}
\begin{split}
&\mathrm{PV}^{-1} \sum_l 
(\rho^i (T_l)- \rho^i (t))(\rho^j (T_l)- \rho^j (t))
\mathbf{CF} (T_l) e^{-(T_l-t) \mathbf{r}_t (T_l) }, \\
& \hspace{6cm} 1 \leq i,j \leq n,
\end{split}
\end{equation}
are generalized convexities 
in the sense of section \ref{sensitivity}.
\end{theorem}
Note that when $ \rho (T) $ is affine in $ T $
$ \iff $ $ \sigma (T) \equiv \sigma $;
i.e., only {\em parallel shifts} are allowed,  
the standard duration and convexity are retrieved. 
\begin{proof}
Using (\ref{interpolation1}) 
we can represent the spot rate $ \mathbf{r} $
by the estimated volatility matrix (and average):
\begin{equation}\label{spotR1}
\begin{split}
\mathbf{r}_t (T_i) 
&= - \frac{1}{T_i- t} \log P^{T_i}_t
= -\frac{\Delta t }{T_i- t} \sum_{u=t}^{Ti -\Delta t} F_t (u) \\
&= - \frac{T_1-t}{T_i- t} F_t (t, T_1) 
- \frac{T_2-T_1 }{T_i-t} F_t (T_1, T_2) - \cdots
- \frac{T_i -T_{i-1} }{T_i-t} F_t (T_{i-1}, T_i ) \\
&= - \frac{ 1 }{T_i- t} \langle\Delta t 
\sum_{u=t}^{T_1 - \Delta t} \boldsymbol{\sigma} (u) 
+ \sum_{l=1}^{i-1}  (T_{l+1} - T_l ) \boldsymbol{\sigma} (T_l, T_{l+1} ) 
,  i_t -i_0 \rangle \\
&=  - \frac{ 1 }{T_i- t} \langle (T_1 -t ) \boldsymbol{\sigma} (T_1, T_2) 
+ \sum_{l=1}^{i-1}  (T_{l+1} - T_l ) \boldsymbol{\sigma} (T_l, T_{l+1} ) 
,  i_t -i_0 \rangle \\
& \hspace{6cm} - \frac{1}{T_i -t} \log P^{T_i}_0 \\
&= - \frac{ 1 }{T_i- t} \{ \langle \boldsymbol{\rho} (T_l) 
- \boldsymbol{\rho} (t)
,  i_t -i_0 \rangle + \log P^{T_i}_0 \}.
\end{split}
\end{equation}
(Here we meant $ \boldsymbol{\sigma} = (\mu,\sigma) $
and $ \boldsymbol{\rho} = (\rho^0, \rho) $. )
The last equality, or 
\begin{equation}\label{duration101}
(T_1 -t) \boldsymbol{\sigma} (T_1, T_2) 
+ \sum_{l=1}^{i-1}  (T_{l+1} - T_l ) \boldsymbol{\sigma} (T_l, T_{l+1} )
= \boldsymbol{\rho} (T_i) - \boldsymbol{\rho} (t),
\end{equation}
is obtained by recalling (\ref{volatility100}). 

Therefore we have
\begin{equation}\label{Ito01}
\begin{split}
& (T_l-t-\Delta t) \mathbf{r}_{t+\Delta t} (T_l)  - (T_l-t) 
\mathbf{r}_{t} (T_l) \\
& \quad = - \langle \boldsymbol{\rho} (T_l) 
- \boldsymbol{\rho} (t), \Delta i_t \rangle 
+ \langle  \boldsymbol{\rho} (t+\Delta t) - \boldsymbol{\rho} (t) ,
i_{t+\Delta t} - i_0 \rangle \\
& \quad = - \langle \boldsymbol{\rho} (T_l) 
- \boldsymbol{\rho} (t), \Delta i_t \rangle 
+ \langle \Delta t \, \boldsymbol{\sigma} (t), 
i_{t+\Delta t} - i_0 \rangle \\
& \quad
= - \langle \rho(T_l) - \rho(t) , 
\Delta \mathbf{w}_t \rangle \\
& \hspace{2cm} + \big\{ 
\langle \boldsymbol{\sigma} (T_1, T_2), 
i_{t+\Delta t} - i_0 \rangle 
- \rho^0(T_l) + \rho^0 (t) \big\} \Delta t \\
& \quad
= - \langle \rho(T_l) - \rho(t) , 
\Delta \mathbf{w}_t \rangle \\
& \hspace{1cm} + \big\{ 
\langle {\sigma} (T_1, T_2), 
\mathbf{w}_{t+\Delta t} - \mathbf{w}_0 \rangle 
- \rho^0(T_l) + \rho^0 (t) + \rho^0 (t+\Delta t)(t+\Delta t) \big\} \Delta t .
\end{split}
\end{equation}
In particular, the left-hand-side of (\ref{Ito01}) 
has the other of $ \sqrt{\Delta t} $.
Therefore,  the asymptotic expansion with respect to $ \Delta t $ of
the present value $ \mathrm{PV} $ of a cash flow $ \mathbf{CF} $ is:
\begin{equation}\label{Taylor101}
\begin{split}
& \frac{\Delta (\mathrm{PV})}{\mathrm{PV}} 
 :=  \frac{\mathrm{PV} ( \mathbf{r} + \Delta 
\mathbf{r}_t ; t + \Delta t ) -
\mathrm{PV} (\mathbf{r}; t) }
{\mathrm{PV}} \\
& = \sum_l \mathbf{CF} (T_l) \left\{ e^{-(T_l-t-\Delta t) 
\mathbf{r}_{t+\Delta t} (T_l)}  - e^{-(T_l-t) 
\mathbf{r}_{t} (T_l) } \right\}/ \mathrm{PV}_t  \\
& \simeq - \sum_l \mathbf{CF} (T_l) e^{-(T_l-t) \mathbf{r}_t (T_l) }\left\{
(T_l-t-\Delta t) \mathbf{r}_{t+\Delta t} (T_l)  - (T_l-t) 
\mathbf{r}_{t} (T_l)  \right\}/ \mathrm{PV}_t \\
& + \frac{1}{2} \sum_l \mathbf{CF} (T_l) e^{-(T_l-t) \mathbf{r}_t (T_l) }
\left\{(T_l-t-\Delta t) \mathbf{r}_{t+\Delta t} (T_l)  - (T_l-t) 
\mathbf{r}_{t} (T_l)  \right\}^2/ \mathrm{PV}_t \\
& \hspace{6cm} + \mathrm{o} (\Delta t ).
\end{split}
\end{equation}

Combining (\ref{Taylor101}), 
(\ref{Ito01}), and (\ref{duration101}), 
we have (\ref{Taylor8}).
To see (\ref{durations-g}) and (\ref{convexities-g})
are generalized durations and convexities respectively,
just compare (\ref{Taylor2}) and (\ref{Taylor8}).
\end{proof}

\subsection{Another framework based on discrete It\^o formula}
Since we have assumed {\em completeness} ($ s=n $) in $(a2')$, 
we can use a discrete It\^o formula (c.f.\cite{A3, A4})
to obtain perfect equality instead of 
the above asymptotic expansion (\ref{Taylor8}). 
\begin{theorem}
Define 
\begin{equation}\label{duraions-g2}
\begin{split}
\tilde{D}_j &:= \frac{e^{\langle \sigma (T_1, T_2), 
\mathbf{w}_{t} - \mathbf{w}_0 \rangle
 - \{ \rho^0 (T_l ) -  \rho^0 (t)\} 
+ \rho^0 (t+\Delta t)
(t+\Delta t)}}{\Delta t} \\
& \hspace{1cm} 
\times \Ex [\Delta w^j e^{\langle \rho(T_l) 
- \rho(t) +\Delta t \sigma(T_1,T_2) , 
\Delta \mathbf{w}_t \rangle}],\,\, j=1,...,n
\end{split}
\end{equation}
and 
\begin{equation}
\begin{split}
\tilde{D}^2 &:= \Ex [e^{\langle \sigma (T_1, T_2), 
\mathbf{w}_{t} - \mathbf{w}_0 \rangle
 - \{ \rho^0 (T_l ) -  \rho^0 (t)\} 
+ \rho^0 (t+\Delta t)
(t+\Delta t) } \\
& \hspace{2cm} \times e^{\langle \rho(T_l) 
- \rho(t) +\Delta t \sigma(T_1,T_2) , 
\Delta \mathbf{w}_t \rangle 
}- 1]/ \Delta t. \\
\end{split}
\end{equation}
Then 
\begin{equation}\label{ItoTaylor7}
\frac{\Delta (\mathrm{PV}_t)}{\mathrm{PV}_t} 
= - \sum_j \tilde{D}_j \Delta w^j_t 
 + \frac{1}{2}\tilde{D}^2 \, \Delta t.
\end{equation}
\end{theorem}

\begin{proof}
Recall 
\begin{equation*}\label{Ito-1}
\begin{split}
& \frac{\Delta (\mathrm{PV}_t)}{\mathrm{PV}_t} 
 = \sum_l \mathbf{CF} (T_l)  e^{-(T_l-t) 
\mathbf{r}_{t} (T_l) } \left\{ e^{-(T_l-t-\Delta t) 
\mathbf{r}_{t+\Delta t} (T_l) +(T_l-t) 
\mathbf{r}_{t} (T_l)  }  - 1 \right\}/ \mathrm{PV}_t, 
\end{split}
\end{equation*}
and 
\begin{equation*}\label{Ito-001}
\begin{split}
& (T_l-t-\Delta t) \mathbf{r}_{t+\Delta t} (T_l)  - (T_l-t) 
\mathbf{r}_{t} (T_l) \\
& \quad
= - \langle \rho(T_l) - \rho(t) , 
\Delta \mathbf{w}_t \rangle \\
& \hspace{1cm} + \big\{ 
\langle {\sigma} (T_1, T_2), 
\mathbf{w}_{t+\Delta t} - \mathbf{w}_0 \rangle 
- \rho^0(T_l) + \rho^0 (t) 
+ \rho^0 (t+\Delta t)(t+\Delta t) \big\} \Delta t.
\end{split}
\end{equation*}
We may regard these random variables 
as a function on $ S $.
By putting the latter to be $F(s)$; i.e.
\begin{equation*}
\begin{split}
& F(s) := \langle \rho(T_l) - \rho(t) + \Delta t \sigma (T_1,T_2), 
\Delta \mathbf{w} (s) \rangle \\
& \hspace{1cm} + \big\{ 
\langle {\sigma} (T_1, T_2), 
\mathbf{w}_{t} - \mathbf{w}_0 \rangle 
- \rho^0(T_l) + \rho^0 (t) 
+ \rho^0 (t+\Delta t)(t+\Delta t) \big\} \Delta t, 
\end{split}
\end{equation*}
the former is represented as
\begin{equation}\label{LofS}
\frac{\Delta (\mathrm{PV}_t)}{\mathrm{PV}_t} 
 = \frac{1}{ \mathrm{PV}_t}
\sum_l \mathbf{CF} (T_l)  e^{-(T_l-t) 
\mathbf{r}_{t} (T_l) } \left\{ e^{F(s)}  - 1 \right\}. 
\end{equation}
Since on the other hand 
$ (1,\Delta w^1/\sqrt{\Delta t},...,\Delta w^n/\sqrt{\Delta t}) $ 
is an orthonormal basis of $ L(S) $, we have
\begin{equation}\label{fourier01}
e^{F(s)} = \sum_{j} \frac{\Ex [\Delta w^j e^F ]}{\Delta t} \Delta w^j
+ \frac{\Ex [e^F]}{\Delta t} \Delta t. 
\end{equation}
Substituting the right-hand-side of (\ref{fourier01}) 
for $ e^{F(s)} $ in (\ref{LofS}), we get the expansion 
(\ref{ItoTaylor7}).
\end{proof}

\section{Passage to IKRS models}\label{PSTL}
\begin{quote}
... So it was very difficult to estimate the parameters 
in the original formulation of the mode. 
That was the motivation for HJM. 
We looked at the continuous time limit of Ho-Lee... (R. Jarrow \cite{Jar})
\end{quote}

\subsection{IKRS interest rate models}
In this section we will show that the forward rates of  
our generalized Ho-Lee model {\em converge}
to those of a special class of HJM which we call
IKRS model, as $ \Delta t $ tends to zero. 

Let us begin with a quick review of
Inui-Kijima-Ritchken-Sankarasubramanian's (IKRS for short)
interest rate model \cite{IK, RS}.
Above all, 
the class is characterized by the separability of 
the volatility structure of instantaneous forward rates
in a continuous-time HJM framework. 

Recall that HJM is based on 
the semi-martingale representation of 
the instantaneous 
forward rate 
\begin{equation*}
F_t (T) 
:= \lim_{\Delta t \downarrow 0} (\Delta t)^{-1} F_t (T, T+\Delta t).
\end{equation*}
In Brownian cases HJM expression is 
\begin{equation}\label{frate}
dF_t (T) = \sum_i v^i_t(T) (dW^i_t -\lambda^i \,dt)
+ \sum_i \left( v^i_t(T) \int_t^T v^i_t(u)\,du \right) \,dt.
\end{equation}
Here $ W^i $'s are mutually independent Wiener processes and
$ \lambda^i $'s are some adapted processes that correspond 
to so-called {\em market prices of risk}. 
IKRS assumes the following separation of volatility:
\begin{equation}\label{IKRSvol}
v^i_t(T) := \eta^i_t 
e^{-\int_t^T \kappa^i (s)\,ds} \quad (i=1,2,...,n),
\end{equation}
where $ \kappa^i $'s are some deterministic functions 
and $ \eta^i_t $'s are some adapted processes.
The favorable property of IKRS lies in the expression of bond prices:
\begin{equation}\label{IKRSbond}
P_t^T =\frac{P_0^T}{P_0^t}
\exp \sum_i \left\{ - \frac{1}{2} |\Lambda^i_t (T)|^2 \phi^i_t 
+ \Lambda^i_t(T)
[F_0(t)-r^i (t)] 
\right\} 
\end{equation}
where 
\begin{equation}\label{IKRSvol2}
\Lambda^i_t(T) = \int_t^T e^{-\int_t^u \kappa_x\,dx} \,du,
\end{equation}
and
\begin{equation}\label{RSSV}
\begin{split}
dr^i_t &= [\partial_t F_0 (t) + \phi^i_t
+ \kappa(t)(F_0(t)-r^i_t )] dt 
+ \eta^i_t (dW^i -\lambda^i dt),\, r^i_0 = 0 \\
d\phi^i_t &= [|\eta^i_t|^2 - 2 \kappa_t
\phi^i_t ]dt, \, \phi^i_0 = 0.
\end{split}
\end{equation}
The ``state variable" $ \{ r^i, \phi^i\}_{i=1,...,n} $ 
can be, if we impose\footnote{We need to be careful about the linear 
growth condition in (\ref{RSSV})
when we use such Markovian implementations. 
For details, see \cite{A2}. }$ \eta^i \equiv \eta (r^i) $,
a solution to a degenerate stochastic differential equation 
under the risk neutral measure, whose density on $ \F_t
:= \sigma (\{W^1_s,...,W^n_s\}; s \leq t) $
with respect to the real world measure is 
$ \E ( \sum_i \int \lambda^i dW^i )_t $.
For details, see e.g. \cite[section 5.3]{Brigo-Mercurio}.

\begin{remark}{\em 
The separability of volatility structure 
is first discussed in \cite{Jam} by
F. Jamshidian, who called the class {\em Quasi Gaussian}.
}
\end{remark}

\subsection{Reparametrizations of IKRS}
While the above parameterization of original IKRS 
is aimed to have Markovian state variables 
which allow reasonable discretizations, 
our parameterization below gives
a \underline{direct discretization scheme by recombining trees}
though it is limited to the Gaussian cases.

Our new insight starts from the following observation:
\begin{lemma}
Let $ \rho^i (t) = \int_0^t e^{-\int_0^s \kappa^i (u)\,du}\,ds $, 
$ \varphi^i(s) = \eta^i (s) e^{\int_0^s \kappa^i(u)\,du} $, 
and assume $ \lambda^i (t) = - \rho^i (t) \varphi^i (t)  $. 
Then
the forward rates (\ref{frate}) with (\ref{IKRSvol}) 
turns into\footnote{Here {\em dot} 
over a function means its derivative. }
\begin{equation}\label{stationary-conti}
F_t(T) = F_0 (T) + \sum_i \left( \dot{\rho}^i (T) \int_0^t \varphi^i (s)\,dW^i_s
+ \dot{\rho}^i (T) \rho^i (T) \int_0^t |\varphi^i (s)|^2 \,ds \right).
\end{equation}
In particular, when $ \varphi^i \equiv 1 $, 
then $ F_t (T) $ for each $ T $ 
is a process  with stationary independent increments. 
\end{lemma}
\begin{proof}
It is direct if one notice that $ \dot{\rho}(T) 
= e^{-\int_0^T \kappa^i (s) \,ds} $ and hence
$ v^i_t = \eta^i_t 
e^{-\int_t^T \kappa^i (s)\,ds} = \varphi^i (t)  \dot{\rho}(T)  $.
\end{proof}

We now consider (\ref{stationary-conti}) itself to be a new model,
meaning that $ \dot{\rho} $ can be negative, 
which is impossible in the original IKRS.

Since the standard forward rate over $ [T, T'] $ is 
$ \frac{1}{T'-T} ( \log P_t^T - \log P^{T'}_t ) =: F_t (T,T') $
in the continuous-time framework is represented by
the instantaneous ones by
\begin{equation*}
F_t (T, T') = \frac{1}{T'-T} \int_T^{T'} F_t (u) \,du .
\end{equation*}
Therefore, for a given $ T, T' $, 
IKRS in the form of (\ref{stationary-conti}) with $ \varphi^i \equiv 1 $ 
gives
\begin{equation}\label{forwardHJM1}
\begin{split}
F_t (T, T') &= F_0 (T,T') 
+ \left\langle \frac{\rho (T') - \rho (T)}{T'-T}, \mathbf{W}_t 
-\mathbf{W}_0 \right\rangle \\
& \qquad + \frac{ |\rho (T')|^2 
- | \rho (T) |^2 }{T'-T}  \,t.
\end{split}
\end{equation}
We may extend IKRS to any right continuous 
$ \rho $ by (\ref{forwardHJM1})
if we do not care about the instantaneous forward rates.

\subsection{A limit theorem}
Comparing (\ref{forwardHJM1}) for $ T'= T_{i+1} $ and $ T=T_i $
with  (\ref{forward401})-(\ref{driftcond2}),
we have the following result.
\begin{theorem}
Let $ t $ be fixed and $ \Delta t = t/N $.
The forward rates $ \{ F_t (T_i,T_{i+1}); i=1,...,k\} $
given by
(\ref{forward401})
of the multi-factor generalization of Ho-Lee
converge weakly as $ N \to \infty $ to 
the corresponding forward rates of 
IKRS given by (\ref{forwardHJM1}) for the same $ \rho $. 
\end{theorem}
\begin{proof}
It suffices to show 
\begin{equation*}
\lim_{\Delta t \to 0}
\frac{2}{\Delta t}
\log \Ex[\exp { \langle  \rho (T), \Delta \mathbf{w} \rangle} ]
= | \rho(T) |^2.
\end{equation*}
But this follows almost directly since we have
\begin{equation*}
\begin{split}
\Ex[\exp { \langle  \rho (T), \Delta \mathbf{w} \rangle} ]
&= 1 + \Ex [ \langle \rho (T), \Delta  \mathbf{w} \rangle
+ \frac{1}{2} \Ex [ \langle \rho (T), \Delta  \mathbf{w} \rangle^2] + 
\mathrm{o} (\Delta t) \\
&= 1+ \frac{1}{2} \langle \rho(T) \otimes \rho(T), 
\Ex [\Delta  \mathbf{w} \otimes \Delta  \mathbf{w}]\rangle 
+ 
\mathrm{o} (\Delta t) \\
&= 1 + \frac{1}{2}\langle \rho(T) \otimes \rho(T), 
\Delta t I_n \rangle + 
\mathrm{o} (\Delta t) \\
&= 1 + \frac{1}{2}|\rho(T)|^2 \Delta t + \mathrm{o} (\Delta t).
\end{split}
\end{equation*}
Here the expectation of a vector is taken coordinate-wisely,
and $ I_n $ is the $ n \times n $ unit matrix.
As a matter of course, we have used (\ref{cov}).
\end{proof}

\begin{remark}{\em
In this continuous-time framework, 
the duration-convexity formula like (\ref{Taylor8}) 
can be obtained directly from It\^o's formula.  
In fact,  since
\begin{equation}\label{spotR2}
\begin{split}
-(T_i- t)\mathbf{r}_t (T_i) &= \langle {\rho} (T_i) 
- {\rho} (t)
,  \mathbf{W}_t - \mathbf{W}_0 \rangle + \{ |\rho(T_i)|^2 - |\rho(t)|^2 \} t 
 + \log P^{T_i}_0 \\
 &= \sum_j \int_0^t \{ \rho^j (T_l) 
- \rho^j (s) \}\, d{W}^j_s 
+ \int_0^t \{ |\rho(T_i)|^2 - |\rho(s)|^2 \}\,ds  \\
& \hspace{4cm} - \sum_j \int_0^t  W^j \,d \rho_s 
- \int_0^t s \, d \mu_s, 
\end{split}
\end{equation}
we have the following equality by applying It\^o's formula.
\begin{equation*}\label{TaylorIto8}
\begin{split}
\frac{d (\mathrm{PV})_t}{\mathrm{PV}_t} 
& = \sum_j \frac{1}{\mathrm{PV}_t}
\sum_l ( \rho^j (T_l) - \rho^j (t) ) 
\mathbf{CF} (T_l) e^{-(T_l-t) \mathbf{r}_t (T_l) }\, dW^j \\
& + \frac{1}{2} \sum_{i,j} \frac{1}{\mathrm{PV}_t} \sum_{l} 
\{ \rho^i (T_l) - \rho^i (t) \} \{ \rho^j (T_l) - \rho^j (t) \} 
\mathbf{CF} (T_l) e^{-(T_l-t) \mathbf{r}_t (T_l) }\,dt \\
&+\frac{1}{\mathrm{PV}_t} 
\sum_l \big(  \rho^0 (T_l ) - \rho^0 (t) - t \dot{\rho}^0 (t) 
- \sum_j  (W^j - W^j_0) \dot{\rho}^j (t) \big) \cdot \\
& \hspace{4cm} 
\cdot\mathbf{CF} (T_l) e^{-(T_l-t) \mathbf{r}_t (T_l) } \,dt.
\end{split}
\end{equation*}
Thus our generalized durations (\ref{durations-g})
and convexities  (\ref{convexities-g}) survive
the continuous-time framework. 
}
\end{remark}

\section{Concluding remark}
For all refinements cited above, 
the problem of absurd negative interest rates 
remains unresolved.
To overcome the puzzle, we need to further generalize 
Ho-Lee model; from linear to non-linear. 
The results will be published as PART II.

\appendix\section{Appendix}

\subsection{Principal component analysis: a review}\label{PCA}

Let 
$ \mathbf{y}_p \equiv ( {y}^1_p,...,{y}_p^k ) $,
$ p = 1,2,...,N $ be sample data of 
an $ \R^k $-valued random variable $ \mathbf{y} $.
Define covariance matrix $ C:=( c_{l,m} )_{1\leq l,m \leq k }$
by the sample covariances:
\begin{equation}\label{cov1}
c_{l,m} =\frac{1}{N}\sum_{p=1}^N 
( {y}_p ^l -  \overline{{y}}^l )
( {y}_p^m -  \overline{{y}}^m ),
\end{equation}
where $ \overline{\mathbf{y}}^l  
= \frac{1}{N} \sum_{p=1}^N {y}_p^l $: 
the sample average. 
Since $ C $ is a positive definite matrix, it is
diagonalized to 
$ \Lambda \equiv \mathrm{diag} [\lambda_1,...,\lambda_k] $
by an orthogonal matrix 
$ U \equiv [\mathbf{x}_1,...,\mathbf{x}_k ] $;
we have $ U C U^* = \Lambda $.
Here $ \mathbf{x}_j \equiv ( x_{j,1},...,x_{j,k}) $ 
is an eigenvector of the eigenvalue $ \lambda_j $.

We can assume without loss of generality that 
$ \lambda_1 \geq \lambda_2 \geq \cdots \geq \lambda_k \geq 0 $, 
and by ignoring small $ \lambda $'s 
we can choose $ \lambda_1,...,\lambda_n $ 
(meaning that $ \lambda_{n+1},....,\lambda_k $ are set to be zero )
and corresponding 
principal component vectors $ \mathbf{x}_1,...,\mathbf{x}_n $.
so that we have an expression of
\begin{equation}\label{express1}
\mathbf{y} =  \sum_{j=1}^n \mathbf{x}_j \sqrt{\lambda_j} w^j 
+ \overline{\mathbf{y}},
\end{equation}
where $ [w^1,...,w^n]=:\mathbf{w} $ is an $ \R^n $-valued random variable 
with 
\begin{equation}\label{orthogonals}
\Ex [ w^j ] = 0 \quad
\mbox{and}
\quad 
\mathrm{Cov} ( \Delta w^i , \Delta w^j) =
\begin{cases}
1 & i=j \\
0 & i \ne j.
\end{cases}
\end{equation}
By redefining
$ \Lambda^{1/2} 
:= \mathrm{diag} [\lambda^{1/2}_1,...,\lambda^{1/2}_n] $
as $ n \times n $ matrix and 
$ X := [\mathbf{x}_1,...,\mathbf{x}_n] $ as 
$ k \times n $ matrix, we can rewrite (\ref{express1}) as
\begin{equation}\label{express2}
\mathbf{y} = X \Lambda^{1/2} \mathbf{w} 
+ \overline{\mathbf{y}},
\end{equation}

The expression (\ref{express1})/(\ref{express2}) is universal
in the following sense:
If $ \mathbf{y} $ has the expression
of
\begin{equation*}
\mathbf{y} = A \mathbf{\tilde{w}}+ \overline{\mathbf{y}},
\end{equation*}
and if 
$ \mathbf{\tilde{w}} \equiv [\tilde{w}^1,...,\tilde{w}^n ] $ 
has the first two 
moments as (\ref{orthogonals}), then 
$ k \times n $ matrix $ A $
should satisfy $ A = X \Lambda^{1/2} T $
for some orthogonal matrix $ T $.
In this case $ T \mathbf{\tilde{w}} $
again has the moments of (\ref{orthogonals}),
and so, the expression (\ref{express2}) is retrieved.

To be more precise, the above is rephrased as follows;
{\em if the covariance matrix $ C $ of
$ \mathbf{y} -\overline{\mathbf{y}} $
is expressed as $ A A^* $, 
then $ A = X \Lambda^{1/2} T $ 
for some orthogonal matrix $ T $.}
Equivalently, 
$ \Lambda^{-1/2} X^* A $ is an orthogonal matrix.
The last statement can be easily checked if one notice 
that $ X \Lambda^{1/2} \Lambda^{1/2} X^* = C $ 
and that $ X^* X $ is the $ n \times n $ unit matrix.

\subsection{No-arbitrage restriction}\label{NAR}
The no-arbitrage principle in our framework is, 
by definition, the impossibility of
\begin{equation}\label{NA1}
( - \mathrm{PV}_t, \mathrm{PV}_{t+\Delta t} ) > 0 \quad 
\mbox{almost surely}
\end{equation}
for any cash flow $ \mathbf{CF} $. 

Since we have assumed that 
{\bf there are only finite possibilities} in $ (a2) $ or $ (a2') $,
the inequality can be regarded as the one in a finite dimensional
(in fact it is $ s+1 $ -dimensional) 
vector space (in the sense that all the coordinates are positive).  
Therefore, by a separation theorem 
from the finite-dimensional convex analysis (see e.g. \cite{duffie:01, Rock}), 
we know that 
the impossibility of (\ref{NA1}) is equivalent to 
the existence of a strictly positive random variable
$ D_t $ 
such that 
\begin{equation}\label{NA3}
\mathrm{PV}_t (s_{t-1},...,s_1)  
= \sum_{s_t \in S }  D_{t + \Delta t}(s_t; s_{t-1},...,s_1 ) 
\mathrm{PV}_{t+ \Delta t}(s; s_{t-1},...,s_1),
\end{equation}
for any $ (s_{t-1},...,s_1) $.

We can reduce the relation (\ref{NA3}) in the following way.
\begin{equation*}\label{NA2}
\begin{split}
 & \sum_{T > t}  \mathbf{CF} (T) e^{-\mathbf{r}_t (\cdot) (T_l) (T_l -t)}
= \sum_{s \in S} D_{t + \Delta t}(s; \cdot) 
\sum_{T >t} \mathbf{CF} (T) 
e^{-\mathbf{r}_{t + \Delta t} (T) (s; \cdot)  (T_l -(t+ \Delta t))}\\
& \hspace{5cm} \mbox{for any flow $ \mathbf{CF} $}, 
\end{split}
\end{equation*}
\begin{equation}\label{NA4}
\begin{split}
&\hspace{1cm} \iff  e^{-\mathbf{r}_t (\cdot) (T) (T -t)}
= \sum_{s \in S} D_{t + \Delta t}(s; \cdot ) 
e^{-\mathbf{r}_{t+ \Delta t} (T) (s; \cdot ) (T -(t+ \Delta t))} \\
&\hspace{1cm} \iff  P^{T}_t (\cdot) 
= \sum_{s \in S} D_{t + \Delta t}(s; \cdot ) P^{T}_{t+ \Delta t} (s; \cdot ) 
\end{split}
\end{equation}
By introducing
\begin{equation*}
\delta_t (s;\cdot ) := D_t (s; \cdot ) /
\Pr (\Delta \mathbf{w}_t= s | \F_{t - \Delta t} )(\cdot ), 
\end{equation*}
we can rewrite the last equation in (\ref{NA4}) as
\begin{equation}\label{centralexp1}
P^T_{t} = \Ex [ \delta_{t+\Delta t} P^T_{t+\Delta t} | \F_t ],
\end{equation}
where $ \F_t = \sigma ( s_t, s_{t-1},...,s_1) $.
By the tower property of conditional expectation 
we have
\begin{equation}\label{NA5}
P^T_{t} = \Ex [ \prod_{u=t+\Delta t}^T \delta_{u}  | \F_t ],
\end{equation}
and in particular 
\begin{equation*}
P^{t+\Delta t}_t = \Ex [ \delta_{t+\Delta t} | \F_t ].
\end{equation*}
By introducing 
\begin{equation*}
\begin{split}
\pi_t (s; \cdot ) & 
:= \frac{\delta_t }{ \Ex [ \delta_{t+\Delta t} | \F_t ]} \Pr 
( \Delta \mathbf{w}_{t+\Delta t} = s | \F_t) (\cdot) \\
&= \frac{\delta_t }{ P^{t+\Delta t}_t} \Pr 
( \Delta \mathbf{w}_{t+\Delta t} = s | \F_t) (\cdot)
= D_t (s; \cdot )/ P^{t+\Delta t}_t, 
\end{split}
\end{equation*}
we have 
\begin{equation}\label{NA6}
P_t^T =  \sum_{s \in S} \pi_{t + \Delta t} (s) P^{t+\Delta t}_t 
P^T_{t + \Delta t}. 
\end{equation}
This means that $ \pi $ is a risk neutral probability 
in the standard terminology: 
Defining $ H_{t+\Delta t} (T) 
:= P^{(t+\Delta t)}(T) P^{(t)} (1)/ P^{(t)} (T+1) $, 
(\ref{NA6}) is equivalent to 
\begin{equation}\label{NA7}
\Ex^\pi  [ H_{t+\Delta t} | \F_t ] = 1,
\end{equation}
where $ \Ex^\pi $ stands for expectation with respect to $ \pi $.

It is easy to see that 
the converses, starting from existence of $ \pi $,  ending in 
that of $ D $, are true.

To summarize, we have the following.
\begin{theorem}\label{NAS}
The following statements are equivalent.
\begin{enumerate}[(i)]
\item The market is arbitrage-free.
\item There exists a strictly positive adapted process $ \{ D_t \}$
such that (\ref{NA4}) holds.
\item There exists a strictly positive adapted process $ \{ \delta_t \}$
such that (\ref{NA5}) holds.
\item There exists a risk neutral probability $ \pi $ 
(such that 
(\ref{NA6}) or (\ref{NA7}) holds). 
\end{enumerate}
\end{theorem}

Now we give a proof of Theorem \ref{SPD}.
\begin{proof}[Proof of Theorem \ref{SPD}]
Suppose we have a state price density $ \mathsf{D} $.
Define $ \delta_t := D_t/D_{t-1} $. Then (\ref{NA5}) holds.
The converse is trivial. 
\end{proof}


\begin{thebibliography}{00}
\bibitem
{A2}
Akahori, J.:
Explosion tests for stochastic integral equations 
related to interest rate models,
{\em J. Math. Sci. Univ. Tokyo} {\bf 5} (1998) 727--745. 

\bibitem
{A3}
Akahori, J.:
A discrete Ito calculus approach to 
He's framework for multi-factor discrete market, 
preprint, Ritsumeikan University
, 2004.

\bibitem
{A4}
Akahori, J.:
Discrete Ito Formulas and Their Applications to Stochastic Numerics,
{\em RIMS kokyuroku} {\bf 1462} 
202-210 (2006) ; proceedings of the 7th Workshop on Stochastic Numerics; Jun 27--29, 2005, RIMS, Kyoto. arXiv math.PR/0603341

\bibitem
{BGM}
Brace, A., Gatarek, D. and Musiela, M.:
The market model of interest rate dynamics, 
{\em Mathematical Finance} {\bf 7} (2) (1997) 127-155.

\bibitem
{Brigo-Mercurio}
Brigo, D., and Mercurio, F.:
{\em Interest Rate Models---Theory and Practice}. 
Springer Finance. Springer-Verlag, Berlin, 2001

\bibitem
{CLM}
Campbell, J.Y., Lo, A.W. and MacKinlay, A.C.:
{\em The Econometrics of Financial Markets},
Princeton University Press, 1997.

\bibitem{CFP}
Chen, L., D. Filipovi\'c, and H. V. Poor: 
Quadratic term structure
models for risk-free and defaultable rates,
{\bf Math. Finance}, 14(4), (2004) 515-536.


\bibitem{CIR}
Cox, J.C., J.E. Ingersoll, and S.A. Ross:
A Theory of the Term Structure of Interest Rates,
{\em Econometrica} {\bf 53} (1985) 385-408.


\bibitem
{duffie:01}
Duffie, D.:
{\it Dynamic Asset Pricing Theory}, 
3rd eds. Princeton University Press, 2001.


\bibitem{D-K}
Duffie, D. and Kan, R.:
A Yield Factor Model of Interest Rates, 
{\em Mathematical Finance} {\bf 6}(4), (1996) 379-406.



\bibitem
{HJM1}
Heath, D., R.A. Jarrow and A. Morton:
Bond Pricing and the Term Structure of Interest Rates: A Discrete Time Approximation. {\em Journal of Financial Quantitative Analysis} 
{\bf 25} (1990), 419-440. 

\bibitem
{HJM2}
Heath, D., R.A. Jarrow and A. Morton:
Bond Pricing and the Term Structure of Interest Rates: 
A New Methodology. {\em Econometrica} {\bf 60} (1992), 77-105.

\bibitem
{HL}
Ho, T.S.Y., and S. B. Lee:
Term-structure movements and pricing interest rate contingent claims,
{\em Journal of Finance} {\bf 41} (1986) 1011--1029.

\bibitem{HL2}
Ho, T.S.Y., and S. B. Lee:
A Closed-Form Multifactor Binomial Interest Rate Model
Ho, Thomas S.Y. Lee, Sang Bin
{\em The Journal of Fixed Income}
{\bf 14-1}, (2004) 8--16.

\bibitem
{Hull}
Hull, J.C.:
{\em Options, Futures, and Other Derivatives}. 6th ed.
Prentice-Hall International Editions. Upper Saddle River, 
NJ: Prentice Hall.
2003.

\bibitem{H-W}
Hull, J. and A. White:
Pricing Interest-Rate Derivative Securities,
{\em Review of Financial Studies}, 
{\bf 3} (1990), 573--592.

\bibitem
{IK}
Inui, K. and Kijima, M.:
A Markovian framework in multi-factor Heath-Jarrow-Morton models,
{\em Journal of Financial and Quantitative Analysis}, {\bf 33} (3), 
(1998), 423 -440.

\bibitem
{Jam}
Jamshidian, F.:
Bond and option evaluation in the Gaussian interest rate model,
{\em Research in Finance}, {\bf 9} (1991), 131--170.

\bibitem
{Jar}
Jarrow, R. interview in 
{\em Derivatives Strategy}, April, 1997.
(available at \url{http://www.derivativesstrategy.com/magazine/archive/1997/0497qa.asp})

\bibitem{Li-Sch}
Littermanm, R. and J. Scheinkman:
Common factors affecting bond returns
{\em The Journal of Fixed Income} {\bf 1} (1988), 54--61.

\bibitem{Lo-Sch}
Longstaff, F and E. S. Schartz:
Interest Rate Volatility and the Term Structure:
A Two-Factor General Equilibrium Model, 
{\em The Journal of Finance} 
{\bf 47} (1992) 1259-1282. 


\bibitem
{Rock}
Rockefeller, RT: 
{\em Convex Analysis}. Princeton University Press, 1970.


\bibitem
{Ross}
Ross, S.A.:
The arbitrage theory of capital asset pricing, 
{\em Journal of Economic Theory}, {\bf 13} (1976) 341-360.

\bibitem
{RS}
Ritchken P. and L. Sankarasubramanian:
Volatility structure of forward rates 
and the dynamics of the term structure,
{\em Mathematical Finance} {\bf 5} (1995) 55--72.


\bibitem{Vas}
Vasicek, O. :
An equilibrium characterisation of the term structure.
{\em Journal of Financial Economics} {\bf 5} (1977) 177--188. 


\end{thebibliography}
\end{document}